\documentclass[a4paper,11pt]{amsart}

\usepackage{amsmath}
\usepackage{amssymb}
\usepackage{amstext}
\usepackage{amsbsy}
\usepackage{amsopn}
\usepackage{amsthm}
\usepackage{amscd}
\usepackage{amsxtra}
\usepackage{amsfonts}

\usepackage{ifthen}
\usepackage{xspace}
\usepackage{fancyheadings}
\usepackage{layout}
\usepackage{hyperref}

\newcommand{\preprintserver}[2]{\href{http://xxx.lanl.gov/abs/math/#2}{#1/#2}}

\input xy
\xyoption{matrix}
\xyoption{arrow}
\xyoption{curve}
\newdir{ (}{{}*!/-5pt/@^{(}}
\newcommand{\xycenter}[1]{\begin{center}
                          \mbox{\xymatrix{#1}}
                          \end{center}
                         }

\newboolean{xlabels}
\newcommand{\xlabel}[1]{
                        \label{#1}
                        \ifthenelse{\boolean{xlabels}}
                                   {\marginpar{#1}}
                                   {}
                       }

\newcommand{\CZ}{\mathbb{C}}

\newcommand{\PZ}{\mathbb{P}}

\newcommand{\HZ}{\mathbb{H}}

\newcommand{\Pn}{\mbox{$\mathbb{P}^n$}}

\newcommand{\sE}{{\mathcal E}}
\newcommand{\sF}{{\mathcal F}}
\newcommand{\sG}{{\mathcal G}}
\newcommand{\sH}{{\mathcal H}}

\newcommand{\sL}{{\mathcal L}}

\newcommand{\sN}{{\mathcal N}}
\newcommand{\sO}{{\mathcal O}}
\newcommand{\sP}{{\mathcal P}}

\newcommand{\sT}{{\mathcal T}}

\newcommand{\sHom}{{\mathcal Hom}}

\newcommand{\PE}{{\PZ(\sE)}}

\newcommand{\suchthat}{\, | \,}

\newboolean{probleme}
\newcommand{\problem}[1]
           {\ifthenelse{\boolean{probleme}}
                       {{\bf(PROBLEM: #1)\bf}}
                       {}
           }
\newboolean{zukuenftiges}
\newcommand{\zukunft}[1]
           {\ifthenelse{\boolean{zukuenftiges}}
                       {{\bf(AUSBAUM\"OGLICHKEIT: #1)\bf}}
                       {}
           }
\newboolean{extras}
\newcommand{\extra}[1]
           {\ifthenelse{\boolean{extras}}
                       {{\bf EXTRA #1 EXTRA\bf}}
                       {}
           }

\DeclareMathOperator{\tr}{tr}

\DeclareMathOperator{\Jac}{Jac}

\DeclareMathOperator{\Tor}{Tor}
\DeclareMathOperator{\rank}{rank}

\DeclareMathOperator{\Sec}{Sec}

\DeclareMathOperator{\spann}{span}
\DeclareMathOperator{\ch}{ch}

\theoremstyle{plain}
\begingroup

\newtheorem{thm}{Theorem}
\newtheorem{cor}[thm]{Corollary}
\newtheorem{lem}[thm]{Lemma}

\newtheorem{prop}[thm]{Proposition}
\newtheorem{conj}[thm]{Conjecture}

\newtheorem{question}[thm]{Question}
\numberwithin{thm}{subsection} 
\endgroup

\newtheorem*{thm*}{Theorem}
\newtheorem*{conj*}{Conjecture}
\newtheorem*{verm*}{Vermutung}

\theoremstyle{definition}
\newtheorem{defn}[thm]{Definition}
\newtheorem{rem}[thm]{Remark}
\newtheorem{example}[thm]{Example}

\numberwithin{equation}{section}

\newcommand{\nosubsections}{\renewcommand{\thethm}{\thesection.\arabic{thm}}
                            \setcounter{thm}{0}
                           }

\newcommand{\cref}[3]{(\ref{#1}, #2 \ref{#3})}

\usepackage{amssymb,amsmath}
\setboolean{xlabels}{false}
\setboolean{probleme}{false}

\begin{document}

\title[Elliptic Normal Curves and their secant varieties]
       {Geometric Syzygies of Elliptic Normal Curves and
       their secant varieties}

\author{Hans-Christian Graf v. Bothmer \and Klaus Hulek}
\address{Institut f\"ur Mathematik\\ Universit\"at Hannover\\ Postfach 6009\\ 30060 Hannover}
\email{bothmer@math.uni-hannover.de, hulek@math.uni-hannover.de}
  
\maketitle

\begin{abstract}
We show that the linear syzygy spaces of elliptic normal curves, their secant varieties and
of bielliptic canonical curves are spanned by geometric syzygies.
\end{abstract}

\newcommand{\spans}[1]{\spann(#1)}

\newcommand{\Secd}{{\Sec_d}}
\newcommand{\Secdminus}{{\Sec_{d-1}}}
\newcommand{\tildeSecdminus}{\widetilde{\Secdminus}}
\newcommand{\sOE}{\sO_E}
\newcommand{\Jacd}{{\Jac_d}}
\newcommand{\sOJacd}{\sO_\Jacd}
\newcommand{\sOEJacd}{\sO_{E \times \Jacd}}
\newcommand{\sOSecd}{\sO_\Secd}
\newcommand{\sOSecdminus}{\sO_\Secdminus}

\nosubsections
\section{Introduction}
Minimal free resolutions of projective varieties 
$X \subset \Pn$ have received considerable
attention in the last years. The aim is to
understand the connection between the geometry of $X$ and the minimal 
free resolution of its ideal $I_X$. In 1984 Green and Lazarsfeld \cite{GL1}
found a way to construct some so called {\sl geometric syzygies} of $I_X$ 
from certain linebundles 
$\sL$ on $X$. For canonical curves Green has conjectured, that these 
geometric syzygies 
determine the shape of the minimal free resolution.

\begin{conj}[Green]
Let $C \subset \PZ^{g-1}$ be a canonical curve without geometric $p$th 
linear syzygies.
Then $C$ as no $p$-th linear syzygies at all.
\end{conj}

This was recently proved by Voisin \cite{VoisinK3}, \cite{VoisinOdd}  
for general $k$-gonal curves of genus $g$ with $k \geq g/3$.
and before that by Teixidor \cite{TeixidorGreen} 
for general $k$-gonal curves
with  $k\le (g+7)/3$. We are, therefore, now at a point where it makes
sense  to ask 
more detailed 
questions about
the connection between geometric syzygies and the minimal free resolution. 
A simple dimension count from Brill-Noether theory 
shows that one cannot expect all syzygies to be
geometric. However, one can ask whether the geometric syzygies span, i.e.

\begin{question}
Do the geometric $p$-th syzygies span the 
spaces of all $p$-th syzygies of a given
variety $X \subset \PZ^n$ ?
\end{question}

The answer to this question is known to be positive in some 
cases. For the  $0$-th syzygies (quadrics) of canonical 
curves this is a theorem of
Green \cite{GreenQuadrics}. It is also true for the $1$-st 
syzygies of general curves of genus 
$g\ge 9$ by \cite{HC} and for $2$-nd syzygies of 
general canonical curves of genus $8$ by \cite{HCandRanestad}.
Furthermore Eusen \cite{Eusen} proves this for all syzygies of  
trigonal canonical curves. 
Since one can easily check it for canonical curves which are 
isomorphic to a plane quintic
\cite{HC} the case of Clifford index $1$ is, therefore, solved.

Going in a different direction the answer is also positive
for rational normal curves and rational normal
scrolls (almost by definition). On the other hand it is not true for 
K3-surfaces $S \subset \PZ^g$
with Picard number $1$ since these do not lie on rank $4$ quadrics.

In this paper we give new results in three different directions. 
First we give a positive
answer to the above question for all syzygies of elliptic normal curves.
Secondly we prove at the same time an 
analogous result for
the higher secant varieties of elliptic normal curves, providing further 
examples of higher dimensional
varieties that do have syzygy spaces generated by geometric syzygies. 
Thirdly, we
use our results to solve the question for bielliptic 
canonical curves of Clifford index $2$. 

It is also interesting to look at the 
variety of geometric $p$-th syzygies inside
the projectivised space of all $p$-th syzygies. 
We show that for elliptic normal curves and
their secant varieties these varieties of geometric syzygies 
contain nondegenerate elliptic
scrolls. One might hope  that these are projectively normal, but 
we show that this hardly ever the case.

Our paper is organised as follows. In sections \ref{s-notation},
\ref{s-linearstrands} 
and
\ref{s-geometricsyzygies} we recall the definition and basic 
properties of graded Betti numbers, 
linear strands and geometric syzygies. 
In sections \ref{s-families}, \ref{s-union} and
\ref{s-intersection} we 
consider families $X$  over a base $B$ whose fibres $X_b$ have minimal degree linear strands 
with identical Betti numbers. We introduce the notion of a family of linear strands  in this 
situation and consider the union 
$\cup X= \cup X_b$ and the intersection $\cap X= \cap X_b$ of fibres. Their 
linear strands are closely related to the cohomology of the family of linear strands.
These two constructions are crucial ingredients in our proofs.
In section \ref{s-mfr} we turn to elliptic normal curves and their higher
secant varieties. 
We investigate the geometry of these secant varieties and compute
their minimal free resolution. 
In 
Section \ref{s-poincare} we use this information to construct a 
family of linear strands whose
union of fibres is the $d$-th secant 
variety $\Secd=\Secd E$ of an elliptic curve
$E$ and whose intersection of fibres is $\Secdminus= \Secdminus E$. 
This and
the symmetry of the minimal free resolution of $\Secdminus$ allows us to 
prove the
geometric syzygy conjecture for these secant varieties in Section \ref{s-geosyz}. We also note
there, that the varieties of geometric syzygies considered in the proof 
are almost never
projetively normal. 
This failure of projective normality is explained in Section \ref{s-kernel}. It
turns out that the missing sections of the geometric syzygy varieties 
of $\Secdminus$ make up the minimal free resolution of $\Secd$. Finally in 
Section Section \ref{s-bielliptic}
we prove the geometric syzygy conjecture for bielliptic canonical curves.

\section{Notation} \xlabel{s-notation}
\nosubsections

\newcommand{\Pspace}{{\PZ^{n-1}}}
Throughout this paper $E \subset \Pspace$ will be an elliptic normal
curve of degree $n$. We denote the origin of $E$ by $A$ and without loss of generality we
can assume that $\sO_E(1) \cong \sO_E(nA)$. The Jacobian of $E$ of degree $d$ complete linear
series will be denoted by $\Jac_d := \Jac_d E$. 
We will also
consider the $d$-secant variety $\Sec_d := \Sec_d E$ of $E$ in $\Pspace$.

More generally, if $X \subset \Pspace$ is any projective variety, we denote
its minimal free resolution by
\[
       F^\bullet_X \to I_X  
       \quad\quad \text{or}\quad\quad
       F^\bullet_X[1] \to \sO \to \sO_X
\]
where we consider $F^\bullet_X$ as a bounded cochain complex
\[
       F^\bullet_X \colon 0 \to \dots \to F^{-2}_X \to F^{-1}_X \to F^0_X
\]
with cohomology concentrated in degree $0$.

For any free cochain complex $F^\bullet$ we write $F^i:= \oplus_j F^i_j \otimes \sO(i-j)$ where the $F^i_j$ vector spaces.
Here and in the rest of the paper  $\sO$ will mean $\sO_\Pspace$, similarly we 
set $\Omega^p := \Omega^p_\Pspace$.

\problem{mention relation with $\Tor$ and Greens $K_{pq}$?}

The dimensions
\[
       \beta_{ij} = \dim F^i_j
\]
are called graded Betti  numbers of $X$. 

\problem{evtl. $\beta_{ij}(I_X)$ definieren?}

Sometimes we will write more shortly
\[ 
  F^i= \oplus_j \sO(i-j)^{\beta_{ij}}
\]
or collect the
graded Betti  numbers $\beta_{ij}$ in a so-called Betti  diagram:
\[
      \begin{array}{cccc|c}
      &       &             &            & \\
      \hline
      & \dots & \beta_{-10} & \beta_{00} & \\
      &       & \beta_{-11} & \beta_{01} & \\
      &       &             & \vdots     & \\  
      & \beta_{ij} &          &            & \\       
      &            &          &            & 
      \end{array}
\]
For better readability we will write a dash (``-'') if $\beta_{ij}=0$.

\begin{example}
The rational normal curve $X \subset \PZ^3$ of degree $3$ has minimal free
resolution
\[
        0 \to \sO(-3)^2 \to \sO(-2)^3 \to \sO \to \sO_X \to 0.
\]
The corresponding Betti  diagram is therefore
\[
    \begin{matrix}
           - & - & 1 \\
           2 & 3 & -
    \end{matrix}
\]
Notice that this notation is dual to the one used by the computer program Macaulay
\cite{M2}. To obtain the diagrams calculated by this program 
one has to take the reflection of our diagrams with respect
to a vertical line. 
\end{example}

\section{Linear Strands} \xlabel{s-linearstrands}
\nosubsections

To study the minimal free resolution $F_X^\bullet \to I_X$ of a 
variety $X$ it is often useful
to linearise its information. We will look at
subcomplexes 
of $F^\bullet_X$ whose differentials are given by 
matrices of linear forms.

\begin{defn}
Let $I_X$ be an ideal sheaf on $\Pspace$, and
\[
       F^\bullet_X \to I_X
\]
a minimal free resolution of $I_X$.  We define the 
complex $F^\bullet_X(d)$ by
\[
       F^{i}_X(d) = F^{i}_d \otimes \sO(i) = \sO(i)^{\beta_{i+d,d}}
\]
with the differentials induced from $F^\bullet_X$.
We call
\[
     F^\bullet_X(d)[-d] \to I_X
\]
the {\sl degree $d$ linear strand} of $I_X$, since the differentials
of $F^\bullet_X(d)$ are given by matrices of linear forms, and the 
map to $I_X$ is defined
by polynomials of degree $d$.
\end{defn}

\begin{defn}
Let $s := s(X)$ be the smallest integer such 
that $H^0(I_X(s)) \not= 0$. Then $s$ is called the {\sl postulation}
of $X$. $F_X^\bullet(s)$ is called the
{\sl minimal degree linear strand} of $X$.
\end{defn}

\begin{example}
Let $F^\bullet_X \to I_X$ be the minimal free resolution of the
rational normal curve in $\PZ^3$. Since $I_X$ is generated by 
quadrics, its minimal degree linear strand is $F^\bullet_X(2)$. Its Betti  diagram is
\[
        \begin{matrix}
          2 & 3 & - & - 
        \end{matrix}
\]
\end{example}

\begin{rem}
 $F_X^\bullet(d)[-d]$ is a subcomplex of $F_X^\bullet$. 
\end{rem}
        
\begin{rem}
Notice that by Hilbert's syzygy theorem, the linear strand $F_X^\bullet(d)[-d]$
is zero in degrees smaller than $-n$.
\end{rem}

We can calculate the linear strand of a scheme $X \subset \Pspace$ by
Koszul cohomology:

\begin{lem} \xlabel{l-koszul}
Let $d \le s(X)$. Then we have
\[ 
  F^{i}_{d} = H^0(I_X \otimes \Omega^{-i-d}(-i))
\]
for the degree $d$ linear strand of $X$. 
\end{lem}

\begin{proof}
By Koszul cohomology \cite{GreenKoszul} the vector space $F^{i}_{d}$ is
the middle cohomology of
\[
           \Lambda^{i-d+1} V  \otimes (I_X)_{d-1} 
           \to
           \Lambda^{i-d} V  \otimes (I_X)_{d} 
           \to
           \Lambda^{i-d-1} V  \otimes (I_X)_{d+1}
\]
where $V=H^0(\sO_\Pspace,\sO(1))$. 
Since $d \le s(X)$ the variety $X$ does not lie on a hypersurface of
degree $d-1$ and the first space vanishes. The kernel of the second map
is easily identified as the homology group above \cite[p. 153]{Ehbauer}.
\end{proof}

\begin{rem}
For $d<s(X)$ this lemma only says that the homology groups above vanish. 
This might seem a trivial observation, but via Proposition \ref{p-cupstrand}
this will provide a crucial vanishing theorem for the vector bundles used
in Section \ref{s-geosyz}.
\end{rem}

An interesting class of varieties whose linear strand is even exact is
given by certain determinantal varieties:
     
\begin{defn}
Let $\phi \colon \sG \to \sH$ be a homomorphism of vector bundles on $\Pspace$. Then
we define the $r$-th degeneracy locus of $\phi$ by
\[
      X_r(\phi) = \{ x \in \Pspace \suchthat \rank \phi(x) \le r \}
\]
equipped with its natural scheme structure.
\end{defn}

\begin{prop} \xlabel{p-determinantal}
Let $G$ and $H$ be vector spaces of dimension $g \ge h$ and
$\phi \colon G\otimes \sO(-1) \to H \otimes \sO$ a map whose  
determinantal locus $X=X_{h-1}(\phi) \subset \Pspace$ is of expected 
dimension. Then the minimal free resolution
\[
                F_X^\bullet \to I_X
\]
is linear, more precisely $F_X^\bullet = F^\bullet(h)[-h]$.
\end{prop}

\begin{proof}
$I_X$ is resolved by the Eagon-Northcott complex.
\end{proof}

The complexity of minimal free resolutions $F^\bullet_X$ tends to 
increase with the codimension of $X$. One approach to understand 
$F^\bullet_X$ is therefore to find varieties $X' \subset \Pspace$ 
of smaller codimension that contain $X$. Often their minimal free 
resolutions contain information about $F^\bullet_X$:

\begin{prop} \xlabel{p-subschemes}
Let $X \subset X' \subset \Pspace$ be subschemes of $\Pspace$. Assume that
$X$ and $X'$ have the same postulation $s(X)=s(X') =: s$. Then there exists
a diagram
\begin{center}
\mbox{
\xymatrix{
             F_{X'}^\bullet(s)[-s] \ar[r] \ar@{ (->}[d] & I_{X'} \ar@{ (->}[d] \\
             F_X^\bullet(s)[-s] \ar[r] & I_{X}
         }
      }
\end{center}
induced by the inclusion $X \hookrightarrow X'$. 
\end{prop}

\begin{proof}
Consider the exact sequence
\[
     0 \to I_{X'} \to I_X \to I_{X'/X} \to 0.
\]
Tensoring with $\Omega^{-i-s}(i)$ defines an inclusion
\[
       H^0(I_{X'} \otimes \Omega^{-i-s}(-i)) \hookrightarrow
       H^0(I_X \otimes \Omega^{-i-s}(-i))
\]
which by Lemma \ref{l-koszul} induces an inclusion
\[
        F^i_{s,X'} \hookrightarrow F^i_{s,X}.
\]
This inclusion commutes with the differentials, since these are
induced by the natural map
\[
        \Omega^{-i-s}(-i)  \to \Omega^{-i-1-s}(-i-1) \otimes H^0(\sO(1))
\]
for both varieties.
\end{proof}

This proposition gives a method to prove the non vanishing of
certain syzygy spaces: 

\begin{example}[Green and Lazarsfeld]
Let $X \subset \PZ^{g-1}$ be a 
canonical curve, and assume that
$|D|$ is a pencil of degree $d$ divisors. Consider the vector spaces
$G=H^0(K-D)$ and $H=H^0(D)$ and the map
\[
        \phi \colon G \otimes \sO_{\PZ^{g-1}}(-1) 
                \to H^* \otimes \sO_{\PZ^{g-1}}
\]
induced by the multiplication of sections. The determinantal variety
$X' = X_1(\phi)$ will contain $X$. It is of expected dimension and
cut out by quadrics. By the propositions above, the
minimal free resolution of $X'$ will be a subcomplex of the 
minimal free resolution of $X$. In particular 
$\beta_{-(g-d-1),2}(X) \not=0$, since $\beta_{-(g-d-1),2}(X') \not=0$.
This proves the ``easy direction'' of Green's conjecture for general
canonical curves.
\end{example}

\section{Geometric Syzygies} \xlabel{s-geometricsyzygies}
\nosubsections

\newcommand{\phiL}{\phi_{\sL}}

We shall now explain the concept of geometric syzygies. Let $X \subset \Pspace$
be a projective variety of postulation $s$ and $\sL$ a line bundle on $X$. This defines
a natural vector bundle homomorphism
\[
    \phiL \colon H^0\bigl(\sO(1) \otimes \sL^{-1}\bigr) \otimes \sO(-1) \to H^0(\sL)^*\otimes \sO
\]
If $u_1,\dots,u_g$ is a basis of $H^0\bigl(\sO(1)\otimes \sL^{-1}\bigr)$ and $v_1, \dots, v_h$
is a basis of $H^0(\sL)$ then $\phiL$ is given by the matrix $(u_i v_j)$. We may assume $g \ge h$, since we can exchange the roles of $\sL$ and $\sO(1) \otimes \sL^{-1}$. Over
$X$ this matrix has rank $1$ and hence $X$ is contained in $X_r(\phiL)$ for
$r \ge1$. Notice that $X_{s-1}(\phiL)$ has also postulation $s$ as it is defined
by $(s \times s)$-minors.

\begin{defn}
An $i$th linear syzygy $f \in F^{-i}_{s,X}$ is called {\sl geometric} if
there exists a line bundle $\sL$ on $X$ such that $f$ is in the image
of the map
\[
       F^{-i}_{s,X'} \to F^{-i}_{s,X}
\]
induced by the inclusion $X \hookrightarrow X' := X_{s-1}(\phiL)$.
\end{defn}

Now consider the special case where $ h:= h^0(\sL)$ is equal to the postulation $s$.
Then $X_{s-1}(\phiL)$ is cut out by the maximal minors. Since $(u_i v_j)$ is
a $1$-generic matrix (see for example \cite{linearsectionsofdeterminantal}) the scheme $X_{s-1}(\phiL)$ has expected codimension \cite[Corollary 3.3]{linearsectionsofdeterminantal}, and hence its minimal free resolution is given by the Eagon-Northcott complex. The geometric
interpretation of $X_{s-1}(\phiL)$ is that set theoretically
\[
    X_{s-1}(\phiL) = \bigcup_{D \in |\sL|} \spans{D} \subset \Pspace
\]
We note at this point that the $1$-genericity of $(u_i v_j)$ also implies that $X_{s-1}(\phiL)$ is smooth outside of $X_{s-2}(\phiL)$ which is at least of codimension $2$ \cite[Corollary 3.3]{linearsectionsofdeterminantal}. As a determinantal variety
of expected codimension $X_{s-1}(\phiL)$ is Cohen-Macaulay and together with the above this
shows that it is normal.

\begin{defn}
A syzygy $f \in F^{-i}_{s,X}$ is called {\sl scrollar}, if there exsists a line bundle $\sL$ on $X$ 
with $h^0(\sL)=s$ such that $f$ is geometric with respect to $\sL$.
 \end{defn}
 
\begin{defn}
We will say that {\sl $X$ satisfies the geometric syzygy conjecture
in step $i$}, if $F^{-i}_{s,X}$ is spanned by geometric syzygies.
\end{defn}

The aim of our paper will be to prove the geometric syzygy conjecture
for elliptic normal curves and their secant varieties. As a corollary
we also obtain the geometric syzygy conjecture for bielliptic canonical
curves.

Since one line bundle produces rarely enough
syzygies to span those of $X$, we will consider
families of line bundles, the determinantal varieties associated to them, and the
corresponding families of linear strands.

\section{Families of linear strands} \xlabel{s-families}
\nosubsections

In this section we will explore the above concepts in a relative setting.
Let 
\begin{center}
\mbox{
\xymatrix{
            X \ar@{ (->}[r] \ar[d] 
            & B \times \Pspace \ar[d]^\sigma \ar[r]^\pi
            & \Pspace
            \\ 
            B \ar@{=}[r]
            & B 
         }
     }
\end{center}
be a family of projective schemes over a smooth, complete and irreducible base $B$. 
We shall assume that the postulation $s(X_b)$ of the fiberes has constant value $s$
for all $b \in B$.

\begin{defn}
A complex $\sF^\bullet(s)$ of the form
\[
   \sF^\bullet(s) \colon \dots \to \sF^{-s-1} \boxtimes \sO(-s-1)
                            \to \sF^{-s} \boxtimes \sO(-s)
\]
with $\sF^i$ vector bundles on $B$ together with a map
\[
        \sF^\bullet(s)[-s] \to I_X
 \]
is called a {\sl family of (minimal degree) linear
strands}, if 
\[
      \sF^\bullet(s)[-s] \otimes \sO_{\PZ_b^{n-1}} \to I_{X_b}
\]
is the minimal degree linear stand of $X_b$ for all $b \in B$.
\end{defn}

If $\sF^\bullet(s)[-s] \to I_X$ is exact one can calculate 
the terms of $\sF^\bullet$ by a relative version of Koszul cohomology:

\begin{lem} \xlabel{l-relativekoszul}
If $\sF^\bullet(s)[-s] \to I_X$ is an exact family of 
linear strands, then 
\[
    \sF^{-i} = \sigma_* (I_X \otimes \pi^* \Omega^{i-s}(i))
\]
and all higher direct images vanish, i.e.
\[
    R^q \sigma_* (I_X \otimes \pi^* \Omega^{i-s}(i)) =0 
\]
for $q\ge1$.
\end{lem}

\begin{proof}
We use hypercohomology with respect to $\sigma$. The second hypercohomology
spectral sequence is
\[
       {}^{II}E_2^{pq} = R^p\sigma_*(H^q  (\sF^\bullet)) 
                         \Rightarrow \HZ_\sigma^{p+q}(\sF^\bullet).
\]    
Since $\sF(s)[-s]$ has cohomology $I_X$ concentrated in degree $0$ we have
\begin{align*}
      {}^{II}E_2^{pq} &= R^p\sigma_*(H^q (\sF^\bullet(s)[-s]\otimes \pi^* \Omega^{i-s}(i))) \\
                                 &= \left\{
                                 \begin{matrix}
                                       R^p\sigma_* (I_X \otimes \pi^* \Omega^{i-s}(i)) & \text{if $q=0$} \\
                                       0                                                                              & \text{otherwise}\\
                                 \end{matrix} 
                                 \right.
\end{align*}
i.e. the spectral sequence degenerates and
\[
   \HZ^p_\sigma := \HZ^p_\sigma(\sF^\bullet(s)[-s] \otimes \pi^*\Omega^{i-s}(i))
                              = R^p\sigma_* (I_X\otimes \pi^* \Omega^{i-s}(i)).
\]           
To prove the identities of the lemma we use the first hypercohomology sequence
\[
       {}^{I}E_1^{pq} = (R^q \sigma_* (\sF^\bullet))_p 
                         \Rightarrow \HZ_\sigma^{p+q}(\sF^\bullet).
\]
We start by calculating the higher direct images of $\sF^\bullet(s)[-s]\otimes \pi^*\Omega^{i-s}(i)$
in step $p$. By Hilbert's syzygy theorem $\sF^\bullet(s)[-s]$ has no non zero terms in steps
$p\le -n$. Similariy by definition it has no non zero terms in steps $p>0$. For $-n-1\le p \le 0$
we have
\begin{align*}
{}^{I}E_1^{pq} 
        &=R^q\sigma_*\bigl[ \sF^\bullet(s)[-s] \otimes \pi^* \Omega^{i-s}(i) \bigr]_p \\
        &=R^q\sigma_* \bigl(\sF^{p-s}(p-s) \otimes \pi^* \Omega^{i-s}(i) \bigr) \\
        &=   \sF^{p-s} \otimes R^q\sigma_* \pi^* \Omega^{i-s}(i-s+p).
\end{align*}        
By Lemma \ref{l-cohomology} in the appendix the second factor is non zero
only for $p=-i+s=-q$. Furthermore $R^{i-s} \sigma_* \pi^* \Omega^{i-s}$ is
trivial.  So
\[
       {}^{I}E_1^{pq} =
                                \left\{
                                 \begin{matrix}
                                       \sF^{-i} & \text{if $p=-i+s=-q$} \\
                                       0           & \text{otherwise.}\\
                                 \end{matrix} 
                                 \right.
\]
So this spectral sequence collapses also and gives
\[
        \HZ^{p+q}_\sigma =
                                 \left\{
                                 \begin{matrix}
                                       \sF^{-i} & \text{if $p+q=0$} \\
                                       0           & \text{otherwise.}\\
                                 \end{matrix} 
                                 \right.
\]
Comparing this with our previous calculation completes the proof.
\end{proof}

\begin{example}
The most important exact families of linear strands are the Eagon-Northcott 
complexes associated to maps of vector bundles
\[
     \phi \colon  \sG \boxtimes \sO(-1) \to \sH \boxtimes \sO
\]
of rank $g \ge h$ such that $X = X_h(\phi)$ is of 
expected dimension on every fiber $X_b$. 
The $X_b$ are determinantal varieties as in proposition 
\ref{p-determinantal}. 
\end{example}

We now want to relate the minimal degree linear strands of a family to those 
of the scheme theoretic intersection $\cap X := \bigcap_{b \in B} X_b$ and
the scheme theoretic union $\cup X := \bigcup_{b \in B} X_b$ of the fibres of
$X$. Suppose that $\cap X$, $\cup X$ and $X_b$ all have the same postulation $s$.
Then by Proposition \ref{p-subschemes} 
the minimal degree linear strand
of $\cup X$ is a subcomplex of all minimal degree linear strands in the family, while
all minimal degree linear strands of the family are subcomplexes of the minimal degree
linear strand of
$\cap X$

\begin{center}
\mbox{
\xymatrix{
             F_{\cup X}^\bullet(s)[-s] \ar[r] \ar@{ (->}[d] 
             &  I_{\cup X} \ar@{ (->}[d]
             \\
             F_{X_b}^\bullet(s)[-s] \ar[r] \ar@{ (->}[d] 
             &  I_{X_b} \ar@{ (->}[d]
             \\
             F_{\cap X}^\bullet(s)[-s] \ar[r]  
             &  I_{\cap X} .
             \\
         }
     }
\end{center}

A natural question is, wether we can construct the minimal degree linear strands 
$F_{\cup_X}^\bullet(s)$ and $F_{\cap X}^\bullet(s)$
from the family $\sF^\bullet(s)$.

\begin{rem}
If $X$ is a family of determinantal varieties obtained from a familiy of line bundles
as in Section \ref{s-geometricsyzygies},
the images of each inclusion 
\[
     F^\bullet_{X_b}(s)[-s] \hookrightarrow F^\bullet_{\cap X}(s)[-s]
\]
will consist of geometric syzygies. 
\end{rem}

\section{The union of fibres} \xlabel{s-union}
\nosubsections

We shall consider the union
\[
           \cup X := \bigcup_{b \in B} X_b
\]
with the natural scheme structure given by
\[
        I_{\cup X} := \pi_* I_X.
\]
Appling the
functor $\pi_*$ to the family of linear strands
\[
       \sF^\bullet(s)[-s] \to I_X
\]
we obtain a complex
\[
       \pi_*\sF^\bullet(s)[-s] \to I_{\cup X}
\]
with
\[
       \pi_*\sF^\bullet(s) \colon \dots  \to H^0(B,\sF^{-s-1})\otimes \sO(-s-1)
                                      \to H^0(B,\sF^{-s})\otimes \sO(-s).
\]
This complex is often the minimal degree linear strand of $\cup X$ :

\begin{prop} \xlabel{p-cupstrand}
Let $\sF^\bullet(s)[-s] \to I_X$ be an exact family of linear strands
with $\sF^{-i} = 0$ for $i>s+n$.  
If $\cup X$ has the same postulation $s$, then its minimal degree linear strand
is
\[
         \pi_* \sF^\bullet(s)[-s] \to I_{\cup X}
\]
otherwise $\pi_* \sF^\bullet(s)$ vanishes.
\end{prop}

\begin{proof}
Since no fiber $X_b$ lies on a hypersurface of degree smaller than
$s$, the same is true for the union $\cup X$.
We can therefore calculate the degree $s$ linear strand of $\cup X$ 
via Lemma \ref{l-koszul}
and \ref{l-relativekoszul}
\begin{align*}
         F_{s,\cup X}^{-i} 
         &= H^0(\Pspace,I_{\cup X} \otimes \Omega^{i-s}(i)) \\
         &= H^0(\Pspace \times B,  I_X \otimes \pi^* \Omega^{i-s} (i)) \\
         &= H^0(B,\sigma_* (  I_X \otimes \pi^* \Omega^{i-s} (i) ) \\
         &= H^0(B,\sF^{-i})
\end{align*}
If the postulation of $\cup X$ is larger than $s$, its degree $s$ linear strand is zero.
\end{proof}

\section{The intersection of fibres} \xlabel{s-intersection}
\nosubsections

The linear strand of the intersection of fibres is much harder to control. For this
we consider functor
\[
     R := R^{\dim B} \pi_* ( - \otimes \sigma^* \omega_B)
\]
where $\omega_B$ is the canonical sheaf on $B$. $R$ is right exact 
since the fibres of $\pi$
have dimension $\dim B$.



%

We start with the following lemma.

\begin{lem}
Let $\sF^{\bullet}(s)$ be a family of (minimal degree) linear strands such that
\[
    \sF^{\bullet}(s)[-s+1] \to \sO_{B \times \Pspace} \to \sO_X \to 0
\]
is an exact sequence. Then there is a natural
isomorphism
\[
   \tr \colon  R \sO_X \cong \sO_{\cap{X}}.
\]
\end{lem}

\begin{proof}
Since 
\[
         \sF^{\bullet}(s)[-s+1] \to \sO_{B \times \Pspace} \to \sO_X \to 0
\]
is exact and restriction to the fibres $X_b$ is right exact, we obtain an exact sequence
\[
              \sF^{-s}_b \otimes \sO(-s) \to \sO_{\PZ^{n-1}_b} \to \sO_{X_b} \to 0.
\]
I.e. $I_{X_b}$ is generated by the image of $\sF^{-s}_b\otimes \sO(-s)$. On the other hand applying $R$ yields the exact sequence
\[
              H^{\dim B} (\sF^{-s} \otimes \omega_B) \otimes \sO(-s) 
              \to H^{\dim B} (\omega_B) \otimes \sO_{\Pspace} \to R \sO_X \to 0.
\]
We shall now compare the two exact sequences. For every point $b \in B$ restriction
defines a map
\[
     H^0(B,(\sF^{-s})^*) \to H^0(B, (\sF^{-s}_b)^*).
\]
Dualizing this gives a homomorphism
\[
    \sF^{-s}_b \to H^{\dim B}(B,\sF^{-s} \otimes  \omega_B) = H^0(R \sF^{-s}).
\]
Similarly we obtain
\[
      \sO_{\PZ^{n-1}_b}  \to H^{\dim B}(\omega_B).
\]
Since the image of $\sF^{-s}_b \otimes \sO(-s)$ generates $H^{\dim B} (\omega_B) \otimes I_{X_b}$, the image of  $H^{\dim B} (\sF^{-s} \otimes \omega_B) \otimes \sO(-s)$ equals
$H^{\dim B} (\omega_B) \otimes I_{\cap{X}}$ and hence the epimorphism
\[
               H^{\dim B} (\omega_B) \otimes \sO_{\Pspace} \to R \sO_X \to 0
\]
induces an isomorphism
 \[
                H^{\dim B} (\omega_B) \otimes \sO_{\cap{X}} \to R \sO_X \to 0.
\]
Since $H^{\dim B} (\omega_B)$ is $1$-dimensional this proves the lemma.
\end{proof}

\begin{rem}
Duality theory defines a trace homomorphism
\[
  \tr \colon  R \sO_X \to \sO_{\cap{X}} \to 0.
\]
In this situation this is (up to possibly a scalar) nothing but the isomorphism
from the above lemma. For this reason we already denoted the isomorphism of the
lemma by $\tr$.
\end{rem}

So if we apply $R$ to
the complex
\[
      \sF^\bullet(s)[-s+1] \to \sO_{B \times \Pspace} \to \sO_X
\]
we obtain a complex
\[
      R \sF^\bullet(s)[-s+1] \to \sO \to \sO_{\cap X}
\]
where
\[
      (R \sF^\bullet)_i = H^{\dim B}(\sF_i \otimes \omega_B) \otimes \sO(-i-d)
                        = H^0(\sF^*_i)^*  \otimes \sO(-i-d)
\]
unfortunately it is not so easy to determine, whether this
is the linear strand of $\cap X$. We shall see in Remark \ref{r-notlinearstrand}
and more detailed in Section \ref{s-kernel} that this is not always the case.

Still we have

\begin{prop} \xlabel{p-psi} 
Assume $s = s(X_b) = s(\cap X)$ and
let $F_{\cap X}^\bullet(s)[-s] \to I_{\cap X}$ be the minimal degree linear strand
of $\cap X$. Then there is a natural map of complexes
\begin{center}
\mbox{
\xymatrix{
            R \sF^\bullet(s)[-s] \ar[r] \ar[d]^\psi & I_{\cap X} \ar@{=}[d] \\ 
                F_{\cap X}^\bullet(s)[-s] \ar[r]    & I_{\cap X} . \\ 
         }
     }
\end{center}
\end{prop}

\begin{proof}
Consider the minimal free resolution
\[
                F_{\cap X}^\bullet  \to I_{\cap X}  
\]
of $\cap X$. Since 
\[
            R \sF^\bullet(s)[-s] \to I_{\cap X}  
\]
is projective, the identity lifts to a map of complexes
\begin{center}
\mbox{
\xymatrix{
            R \sF^\bullet(s)[-s] \ar[r] \ar[d]^\psi & I_{\cap X} \ar@{=}[d] \\ 
                F_{\cap X}^\bullet \ar[r]    & I_{\cap X}  \\ 
         }
     }
\end{center}
which is unique up to homotopy equivalence (see for example \cite[p. 35, Porism 2.2.7]{Weibel}).
For degree reasons the image of $\psi$ lies in the minimal degree linear strand
\[
                F_{\cap X}^\bullet(s)[-s]  \to I_{\cap X} . 
\]
Once again by degree reasons there exists no non trivial homotopy equivalence in
this situation.
\end{proof}      

\begin{rem}
If $X$ is a family of determinantal varieties obtained from a family of line bundles as in
Section \ref{s-geometricsyzygies}, then the image of $\psi$
is the part of the minimal free resolution of $\cap X$ that is 
spanned by geometric syzygies coming from this family. Our
plan is to find a family $X$ such that $\cap X = \Sec_{d-1}$
and $\psi$ surjective.
\end{rem}

Before we can do this, we first construct the minimal free resolution
of $\Sec_{d-1}$ as a mapping cone.

\section{the minimal free resolution of $\Sec_{d-1}$}  \xlabel{s-mfr}
\nosubsections

In this section we construct a minimal free resolution of the
secant variety $\Sec_{d-1}$ of $E$. The main idea is to find a scroll
$X_{|D|}$ that contains $\Sec_{d-1}$ such that the ideal sheaf
$I_{\Sec_{d-1}/X_{|D|}}$ is the dualizing sheaf of $X_{|D|}$.
In this situation a minimal free resolution for $I_{\Sec_{d-1}/X_{|D|}}$
is given by the dual of the minimal free resolution of $X_{|D|}$ which is 
itself just an Eagon-Northcott complex. The minimal free resolution 
of $\sO_{\Sec_{d-1}}$ is then obtained by a mapping cone construction
from these two resolutions. We learned this method from \cite{Sch86}.

Let $\sL=\sO_E(D)$ where $D$ is a divisor on $E$  of degree $d$. From now on 
we shall assume $d \le \frac{n}{2}$. 
The multiplication map 
\[
      \mu \colon H^0(nA-D) \otimes H^0(D) \to H^0(nA)
\]
induces a map of vector bundles
\[
      \phi_\sL \colon H^0(nA-D) \otimes \sO(-1) 
           \to H^0(D)^* \otimes \sO
\]
on $\Pspace$. As in Section \ref{s-geometricsyzygies} this defines the scroll
\[
    X_{|D|} := X_{d-1} (\phi_\sL) = \bigcup_{D' \in |D|} \spans{D'} \subset \Pspace.
\]
On the projective
space $\PZ(H^0(D)^*) \cong \PZ^{d-1}$ the multiplication $\mu$ induces 
an injective
map $\tilde{\mu}$ with cokernel $\sE$:
\[
   0 \to H^0(nA-D) \otimes \sO_{\PZ^{d-1}}(-1) 
     \xrightarrow{\tilde{\mu}} H^0(nA) \otimes \sO_{\PZ^{d-1}}
     \xrightarrow{\rho} \sE
     \to 0
\]
Let $\tau \colon \PE \to \PZ^{d-1}$ be the structure map.
$X_{|D|}$ is the image of the projectivization of $\rho$
\[
  \rho \colon \PE \to \Pspace.
\]
Notice that $\PE$ is the canonical desingularisation of the determinantal variety
$X_{|D|}$ given by
\[
        \PE = \{ (x,y) \in \Pspace \times \PZ^{d-1} \suchthat y \circ \phi_x = 0 \}
\]
where $\phi_x \colon H^0(nA-D) \to H^0(D)^*$ is restriction of $\phi$ to $x$ and $y$
is viewed as a linear form $y \colon H^0(D)^* \to \CZ$.

For future calculations we will
denote by $H$ the pullback of the hyperplane class of $\Pspace$ to $\PE$, i.e. 
$\sO_{\PE}(H)= \sO_{\PE}(1)$. Let $R$ be the preimage 
of a hyperplane of $\PZ^{d-1}$.

\begin{rem}
$X_{|D|}$ contains $\Sec_{d-1}$ since for every set of points 
$p_1,\dots,p_{d-1}$ we have
\[
       H^0(D-p_1-\dots-p_{d-1}) = 1
\]
i.e there exist a unique fiber of $\PE$ whose image contains
the secant through $p_1,\dots,p_{d-1}$. Moreover there is only one such
fiber.
\end{rem}

\newcommand{\omegaPE}{K_\PE}

We want to show, that the preimage
\[
      \tildeSecdminus = \rho^{-1} \Sec_{d-1}
\]
and $\omegaPE$ have the
same class. 
For this we first calculate some intersection numbers. 

\begin{lem} Let $X=X_{l-1}(\phi)$ be a determinantal variety
defined by the maximal minors of a $k \times l$-matrix ($k\ge l$)
of linear forms. If $X$ has the expected codimension $k$ then
\[
     \deg X = {k \choose l-1}.
\]
\end{lem}

\begin{proof}
Since $X$ is a linear section of the generic determinantal variety $M_{l-1}(k,l)$ 
and of expected dimension, both varieties have the same degree. This
degree is calculated in \cite[p. 95]{ACGH} as
\begin{align*}
       \deg(M_{l-1}(k,l)) 
       &= \Pi_{i=0}^{l-(l-1)-1} \frac{(k+i)!i!}{(l-1+i)!(k-(l-1)+i)!} \\
       &= \frac{k!}{(l-1)!(k-l+1)!} \\
       &= { k \choose l-1}.
\end{align*}
\end{proof}

\begin{cor}
On $\PE$ we have the following intersection numbers
\[
       H^{2d-i-2}R^i = { n-d \choose d-i-1 }
\]
\end{cor}

\begin{proof} 
The class $R^i$ is realized by the preimage of a codimension $i$
linear subspace $\PZ^{d-i-1} \subset \PZ^d$. Its image in $X_{|D|}$
is defined by the
maximal minors of the corresponding $(n-d)\times(d-i)$-submatrix
of $H^0(nA-D) \otimes H^0(D) \to H^0(nA)$. By the lemma above
its degree is
\[
      H^{2d-i-2}R^i = \deg \pi^* \PZ^{d-i-1} = { n-d \choose d-i-1 }.
\]
\end{proof}

This allows us to calculate the class and degree of the canonical divisor
on $\PE$:

\begin{prop} \xlabel{p-canonicalclass}
The canonical divisor $\omegaPE$ has class
\[
        \omegaPE= -dH + (n-2d)R
\]
and degree 
\[
          H^{2d-1} \omegaPE= 
        - {n-d+1 \choose d-1}
        - {n-d \choose d-2}.
\]
\end{prop}

\begin{proof}
The relative canonical sheaf of $\PE$ is 
\[
     \omega_{\PE/\PZ^{d-1}} = \tau^*(\Lambda^d \sE)(-dH)
                                 = \sO(-dH+(n-d)R)
\]
by the defining sequence of the vector bundle $\sE$. On the other hand
$\tau^*(\omega_{\PZ^{d-1}})= \sO(-dR)$, and therefore
$\omegaPE = -dH+(n-2d)R$.                                 

The degree of $\omegaPE$ is 
\begin{align*}
      H^{2d-1} \omegaPE 
     &= -dH^{2d-2} + (n-2d)H^{2d-3}R \\
     &= -d{n-d \choose d-1} + (n-2d){ n-d \choose d-2} \\
     &= - {n-d+1 \choose d-1}
        - {n-d \choose d-2}.
\end{align*}

\end{proof}

\begin{prop}\xlabel{p-degreesecantvariety}
Let $d < \frac{n}{2} + 1$ then
the degree of the secant-variety $\Sec_{d-1} E$ is 
\[
   \deg \Sec_{d-1} E=  { n-d \choose d-2 }
                 + { n-d+1 \choose d-1 }.
\]
\end{prop}

\begin{proof}
The codimension of $\Sec_{d-1} E \subset \Pspace$ is $n-2d+2$.
If we choose a general $n-2d+2$-dimensional linear subspace
of $\Pspace$ it intersects $\Sec_{d-1} E$ in $\deg \Sec_{d-1} E$ points.
Now project from this subspace. We obtain a degree $n$ elliptic
curve in $\PZ^{2d+4}$ whose number of $(d-1)$-secant $\PZ^{d-3}$'s is
exactly $\deg \Sec_{d-1} E$. This number $v_{d-2}$ is calculated in
\cite[p. 351]{ACGH}:
\begin{align*}
          v_{d-2} &= \sum_{\alpha=0}^{d-1}
                 (-1)^\alpha { 1+2(d-2)-n \choose \alpha }
                             { 1 \choose d-1-\alpha }\\
              &= (-1)^{d-2}      {2d-3-n \choose d-2 }
                 + (-1)^{d-1}{ 2d-3-n \choose d-1 } \\
              &=   { n-d \choose d-2 }
                 + { n-d+1 \choose d-1 }. \\
\end{align*}
\end{proof}                        



 
\begin{cor} \xlabel{c-anticanonical}
The preimage $\tildeSecdminus$ of the secant variety $\Sec_{d-1}$ is an anticanonical divisor
on $\PE$.
\end{cor}

\begin{proof}
Since $\tildeSecdminus$ is a divisor of $\PE$ it has class $a H + b R$,
for some integers $a$ and $b$ to be determined.

By the Propositions \ref{p-canonicalclass} and \ref{p-degreesecantvariety}  we have
\[
     H^{2d-3} \tildeSecdminus = -H^{2d-3} \omegaPE
\]
and  $\omegaPE = -dH + (n-2d)R$. 
So it is enough to show, that $a =d$:

The image of a fiber of $\PE$ in $\Pspace$ contains $d$ points of $E$ and therefore
$d$ distinct $(d-1)$-secant $(d-2)$-planes. 
Conversely the intersection of any other ${d-1}$-secant $(d-2)$-plane 
with the image of a fiber of $\PE$
which is a $d$-secant $(d-1)$-plane is spanned by their common secant points
(see Lemma \ref {l-normalspan} in the appendix), and this span is contained in one
of the $(d-1)$-secants above. This shows  
\[
       (\tildeSecdminus E) . R^{d-1} = dHR^{d-1}
\]
and consequently $a = d$.
\end{proof}

\begin{lem}
Let $F_{|D|}^\bullet[1] \to \sO \to \sO_{X_{|D|}}$ be the minimal 
free resolution of $X_{|D|}$,
then $F_{|D|}^\bullet = F_{|D|}^\bullet(d)[-d] = EN^\bullet(\phi)[-d]$ where
\[
      EN^{-i}(\phi) = \Lambda^{i} H^0(nA-D) 
                      \otimes \Lambda^d H^0(D) 
                      \otimes S_{i-d} H^0(D) 
                      \otimes \sO(-i)
\]
is the Eagon-Northcott complex associated to $\phi$.
\end{lem}

\begin{proof}
Since $X_{|D|}$ is the locus where $\phi$ drops rank, and since
$X_{|D|}$ is of expected dimension, the Eagon-Northcott complex
is a minimal free resolution.
\end{proof}

\begin{example}
Let $n=6$. Then the Betti  diagram of $\sO_{X_{|D|}}$ is
\[
\begin{matrix}
    \begin{matrix}
          - & - & - & 1 \\ 
          3 & 8 & 6 & -\\
    \end{matrix} 
    &
    \quad\quad\text{and}\quad\quad
    &
    \begin{matrix}
          - & 1   \\
          - & -   \\
          1 & - 
    \end{matrix} 
    \\
    \text{if $\deg D = 2$}
    &&
    \text{if $\deg D = 3$.}
\end{matrix}
\]
\end{example}

Before we can prove the main result of this section, we need a technical lemma.

\begin{lem} \xlabel{l-idealomega}
If $\omega_{X_{|D|}}$ is the dualizing sheaf of $X_{|D|}$ then
\[
        I_{\Sec_{d-1}/X_{|D|}} \cong \omega_{X_{|D|}}
\]
\end{lem}

\begin{proof}
Let $U$ be the regular part of $X_{|D|}$. Since the singular locus of $X_{|D|}$ contains the locus 
where $\rho$ has positive dimensional fibres,
the restriction of $\rho$ to $\rho^{-1}(U)$ defines an isomorphism of $\rho^{-1}(U)$ with $U$. Hence
by the above Corrolary \ref{c-anticanonical} we find that
\[
        I_{\Sec_{d-1} \cap U/U} = \omega_U.
\]
Pushing this isomorphism forward with the inclusion
\[
      i \colon U \to X_{|D|}
\]
gives the result since $X_{|D|}$ is normal as seen in Section \ref{s-geometricsyzygies}
\end{proof}

We have now all the necessary ingredients to prove the main result
of this section:
 
\begin{thm} \xlabel{t-mfr}
Let $F_{\Sec_{d-1}}^\bullet \to I_{\Sec_{d-1}}$ be the minimal free resolution
of the $d-1$-secant variety of $E$, and $D$ any
divisor of degree $d \le n/2$. Then the minimal free resolution of
$\sO_{\Sec_{d-1}}$ is
\[       
     \sO(-n) 
     \to  F_{\Sec_{d-1}}^\bullet(d)[-d+1] 
     \to \sO
     \to \sO_{\Sec_{d-1}}
\]
and there exist an exact sequence
\[
      0 \to F_{|D|}^\bullet(d) 
        \to F_{\Sec_{d-1}}^\bullet(d)
        \to F_{|D|}^\bullet(d)^*
        \to 0
\]
where $F_{|D|}(d)^* := 
       \sHom(F_{|D|}^\bullet(d),\omega_\Pspace)[n]$.
\end{thm}

\begin{proof}
By Lemma \ref{l-idealomega} above we have 
$I_{\Sec_{d-1}/X_{|D|}} = \omega_{X_{|D|}}$. If we apply
$\sHom(-,\omega_{\Pspace})[n]$ to the minimal free resolution
\[
     F_{|D|}^\bullet(d)[-d+1]  \to \sO \to \sO_{X_{|D|}}
\]
we obtain an exact sequence
\[
    \sO(-n) \to F_{|D|}^\bullet(d)^*[-d] \to \omega_{X_{|D|}}
\]
i.e. a free resolution of $\omega_{X_{|D|}}$. 

Now consider the short exact sequence
\[
         0 \to I_{X_{|D|}} \to I_{\Sec_{d-1}} \to I_{\Sec_{d-1}/X_{|D|}} \to 0.
\]   
By the horseshoe lemma
the direct sum of the two minimal free resolutions of $I_{X_{|D|}}$ 
and $I_{\Sec_{d-1}/X_{|D|}}$
fits together to form a (not necessarily minimal) free resolution of
the middle term $I_{\Sec_{d-1}}$. For degree reasons this resolution is
minimal. This proves that the minimal free resolution of $\sO_{\Sec_{d-1}}$
has the shape claimed and that for the linear strands of all three resolutions
in this construction the claimed exact sequence hold.
\end{proof}

\begin{rem}
In the case $d= (n+1)/2$ the variety $\Sec_{d-1}$ is a hypersurface of degree $n$. Here the analog of
Theorem \ref{t-mfr} still holds. The minimal free resolution is given by
\[
      \sO(-n) \to \sO \to \sO_{\Sec_d}. 
\]
\end{rem}

\begin{example}
Let $n=6$. Then we obtain the Betti  diagram of $\sO_{E}=\sO_{\Sec_1 E}$ via
\[
      \begin{matrix}
          - & - & - & - & 1\\
          - & 3 & 8 & 6 & -\\
          - & - & - & - & -\\
      \end{matrix}
      \quad + \quad
      \begin{matrix}
          - & - & - & - & -\\
          - & 6 & 8 & 3 & - \\
          1 & - & - & - & - \\
      \end{matrix}
      \quad = \quad
      \begin{matrix}
          - & - & - & - & 1\\
          - & 9 & 16 & 9 & - \\
          1 & - & - & - & - \\
      \end{matrix}
\]
Similarly we obtain the Betti  diagram of $\sO_{\Sec_2{E}}$ via
 \[
      \begin{matrix}
          - & - & 1 \\
          - & - & - \\
          - & 1 & - \\
          - & - & - \\
          - & - & - \\
      \end{matrix}
      \quad + \quad
      \begin{matrix}
          - & - & - \\
          - & - & - \\
          - & 1 & - \\
          - & - & - \\
          1 & - & - \\
      \end{matrix}
      \quad = \quad
      \begin{matrix}
          - & - & 1 \\
          - & - & - \\
          - & 2 & - \\
          - & - & - \\
          1 & - & - \\
      \end{matrix}
\]
\end{example}

The above theorem allows us to give closed formulae for the graded Betti  
numbers of $\Sec_{d-1}$:

\begin{cor} \xlabel{c-secdbetti}
$\beta_{-id} = \dim(F_{\Sec_{d-1}}^{-i})_d = { n-d \choose i}{i-1 \choose i-d}
                               +{ n-d \choose n-i}{n-i-1 \choose n-i-d}$.
\end{cor}

\begin{proof}
\begin{align*}
  \dim(F_{Sec_{d-1}}^{-i})_d &= \dim EN^{-i}  + \dim EN^{n-i} \\
                             &= { n-d \choose i}{i-1 \choose i-d}
                               +{ n-d \choose n-i}{n-i-1 \choose n-i-d}
\end{align*}
\end{proof}

\begin{cor} \xlabel{c-secdgorenstein}
Let $d \le (n+1)/2$. Then
$\Sec_{d-1}$ is an arithmetically Gorenstein variety. Moreover
the dualizing sheaf of $\Sec_{d-1}$ is trivial
\[
      \omega_{\Sec_{d-1}} = \sO_{\Sec_{d-1}}
\]
\end{cor}

\begin{proof}
By the Auslander-Buchsbaum formula we can read of the depth of the coordinate ring
of $\Sec_{d-1}$ from the minimal free resolution. Since this is the same as the codimension
of $\Sec_{d-1}$ in $\Pspace$ this shows that $\Sec_{d-1}$ is arithmetically Cohen-Macaulay.

By applying $\sHom(-,\omega_{\Pspace})[n]$ to the minimal free resolution
of $\sO_{\Sec_{d-1}}$ we obtain a minimal free resolution of 
$\omega_{\Sec_{d-1}}$. By the Theorem \ref{t-mfr} and the remark following it the minimal free 
resolution of $\Sec_{d-1}$ is symmetric with respect to 
$\sHom(-,\omega_{\Pspace})[n]$, hence $\omega_{\Sec_{d-1}} = \sO_{\Sec_{d-1}}$. This
also shows that $\Sec_{d-1}$ is arithmetically Gorenstein.
\end{proof}

\begin{prop} \xlabel{p-secdminussmooth}
Let $d \le (n+1)/2$, then $\Sec_{d-1}$ is smooth outside of $\Sec_{d-2}$.
\end{prop}

\begin{proof}
\newcommand{\sOLambda}{\sO_{\Lambda_D}}
Let $S^{d-1} E$ be the $(d-1)$st symmetric product of the elliptic curve $E$.
Consider the incidence variety
\[
             I := \{ (D,x) \in S^{d-1} E\times \Pspace \suchthat x \in \spans{D} \}
\]
and its projections
\xycenter{
                      I \ar[r]^\nu \ar[d]^\tau
                      & \Pspace \\
                      S^{d-1} E.
                }
Since each divisor $D \in  S^{d-1} E$ spans a $\PZ^{d-2}$ this incidence variety is
a projective bundle over $S^{d-1} E$ and therefore smooth. The singular locus of $\Secdminus$
is therfore the union of the locus where $\nu$ is not $1:1$ and the image of
the locus in $I$ where the differential $d\nu$ is not injective. 

By Lemma \ref{l-normalspan} the intersection of two $(d-1)$-secant $\PZ^{d-2}$s has to be contained
in a $(d-2)$-secant $\PZ^{d-3}$, so $\nu$ is $1:1$ outside the preimage 
of $\Sec_{d-2}$. So we are left
to consider the differential of $\nu$. For this we look at a divisor $D \in S^{d-1} E$ and the 
$d-1$-secant $\PZ^{d-2} =: \Lambda_D$ spanned by $D$. From now on we identify $\Lambda_D$ with its preimage $\nu^{-1} \Lambda_D$ which is the fiber of $I$ over the point $D \in S^{d-1} E$.
On $\Lambda_D$ the tangent bundle
of $I$ decomposes as
\[
            \sT_I|_{\Lambda_D} =  T_{S^{d-1} E,D} \otimes \sOLambda \oplus \sT_{\Lambda_D}.
\]
Where $T_{S^{d-1} E,D}$ is the tangent {\sl space} to $S^{d-1} E$ in $D$ and $\sT_{\Lambda_D}$
is the tangent {\sl bundle} of the fiber $\Lambda_D$ over $D$.

Since $\nu$ embedds $\Lambda_D$ into $\Pspace$,  we have a diagram
\xycenter{
    0 \ar[r]
    & \sT_{\Lambda_D}  \ar@{ (->}[r] \ar@{=}[d]
    & \sT_I|_{\Lambda_D} \ar[r] \ar[d]^{d\nu}
    & T_{S^{d-1} E,D}  \otimes \sOLambda \ar[r] \ar[d]^{d\nu'}
    & 0
    \\
    0 \ar[r]
    & \sT_{\Lambda_D} \ar@{ (->}[r]
    & \sT_{\Pspace}|_{\Lambda_D} \ar[r]
    & \sN_{\Lambda_D/\Pspace} \ar[r]
    & 0
    }
(A similar construction is used in \cite[section 2]{FocalSystems}).
For us it is now enough to show that $d\nu' $
is injective outside of the preimage $\nu^{-1}(\Sec_{d-2} \cap \Lambda_D)$. 
For this we give an explict description
of $d\nu'$. First recall that $T_{S^{d-1} E,D}$ is canonically isomorphic to $H^0(\sO_D(D))$.
(See for example \cite[p. 160]{ACGH}). Also the linear forms on $\Pspace = \PZ(H^0(\sOE(nA))$
that vanish on $D$ are given by the cokernel of the map
\[
       H^0(\sOE(nA-D)) \hookrightarrow H^0(\sOE(nA))
\]
induced by $D$. This cokernel is easily identified as $H^0(\sO_D(nA))$. Consequently we have
\[
        d\nu'  \colon H^0(\sO_D(D)) \otimes \sOLambda \to H^0(\sOE(nA-D))^* \otimes \sOLambda(1)
\]
and this is induced by the restriction and multiplication of sections
\xycenter{
                     H^0(\sO_D(D)) \otimes H^0(\sOE(nA-D)) \ar[r] \ar[d]
                     & H^0(\sO_D(nA)) \\
                     H^0(\sO_D(D)) \otimes H^0(\sO_D(nA-D)). \ar[ru]
                }
 Since $2(d-1) < n$ the restriction is surjective and $d\nu'$ factors as
 \xycenter{
                     H^0(\sO_D(D)) \otimes \sOLambda \ar[r]^-{d\nu'} \ar[dr]_-{d\nu''}
                      & H^0(\sOE(nA-D))^* \otimes \sOLambda(1) 
                     \\
                     & H^0(\sO_D(nA-D))^* \otimes \sOLambda(1). \ar@{ (->}[u]
                }
To prove injectivity of $d\nu'$ it is enough to prove it for $d\nu''$. This map can be represented
by a square matrix of linear forms. We will choose appropriate bases to do so. 
 
For this write $D = D_1 + \dots + D_m$ with $D_i$ divisors supported on a single point
$P_i$ and decompose 
\[
    H^0(\sO_D(D)) \cong H^0(\sO_D(nA-D))=: V_1 \oplus \dots \oplus V_m
\]
where $V_i$ contains those sections that are supported only on $D_i$. The matrix $M$
representing $d\nu''$ with respect to this decomposition is made up of block matrices $M_{ij}$
representing the multiplication 
of sections
\[
           d\nu''_{ij}  \colon V_i \otimes V_j \to H^0(\sO_D(nA)).
\]
By the choice of $V_i$ this is zero for $i \not= j$, i.e. the block matices outside the diagonal
vanish. Therfore
the determinant $\det M$ is just the product of the determinants $\det M_{ii}$. Lets consider one
such $M_{ii}$. $V_i$ has dimension $v = \deg D_i$ and a basis $\{1,t,\dots,t^{v-1}\}$ of
sections supported on $D_i$. The multiplication of sections $d\nu''_{ii}$ is then represented by
\[ 
        M_{ii} = 
        \begin{pmatrix}
               1          &   t     & \cdots & t^{v-1} \\
               t           & \cdots  & t^{v-1} & 0 \\
               \vdots &  &  & \vdots\\
               t^{v-1} & 0 & \cdots & 0 
        \end{pmatrix}
\]
where the $t^i \in H^0(\sO_D(nA))$ are considered as linear forms on $\Lambda_D$. In particular
$t^{v-1}$ is the linear form that vanishes on $P_i$ with multiplicity $v-1$ and on all points of $P_j$ with the maximal multiplicity $\deg D_j$. This cuts out the corresponding
$(d-2)$-secant in $\Lambda_D$. So $M_{ii}$ has maximal rank outside this $(d-2)$-secant.
This shows that $M$ has maximal rank outside of $\nu^{-1}(\Sec_{d-2} \cap \Lambda_D)$ as claimed.
\end{proof}

\begin{cor} \xlabel{c-secdnormal}
Let $d \le (n+1)/2$, then $\Sec_{d-1}$ is normal.
\end{cor}

\begin{proof}
$\Sec_{d-1}$ is Cohen-Macaulay and smooth in codimension $1$.
\end{proof}

\section{Poincar\'e linebundle} \xlabel{s-poincare}
\nosubsections

Consider the Poincar\'e bundle $\sP$ on $E \times \Jacd$ which is uniquely determined
by the normlization $\sP |_{\Jacd \times {\{A\}}} \cong \sOJacd$. Using the diagram 
\xycenter{
                     &E \times \Jac_d \ar[dl]_p \ar[dr]^q \\
                     E && \Jac_d
}
we define two vector bundles on $\Jacd$ by setting
\[
     \sG := q_* \bigl(p^* \sOE(nA) \otimes \sP^{-1}\bigr)
      \qquad \text{and} \qquad 
      \sH := q_* \sP.
\]
The fiber of $\sG$ over a point $\sOE(D) \in \Jacd$ is $H^0(nA-D)$ and the 
fiber of $\sH$ over this point is $H^0(D)$. Fiberwise multiplication defines a
vector bundle homomorphism
\[
         \sG \otimes \sH \to H^0(nA) \otimes \sOE.
\]
Now projectivize $H^0(nA) \otimes \sO_E$ to obtain a map
\[
    \phi \colon \sG \boxtimes \sO(-1) 
                 \to \sH^*\boxtimes \sO 
\]
on
\begin{center}
\mbox{
\xymatrix{
             \Jacd \times \Pspace \ar[r]^-\sigma \ar[d]^\pi & \Jac_d \\
             \Pspace .
         }
     }
\end{center}
Let $X := X_{d-1}(\phi) \in \Jacd \times \Pspace$ be the submaximal rank locus of $\phi$.
The fiber over a point $\sOE(D) \in \Jacd$ is the scroll $X_{|D|}$ and hence $X$ is of expected dimension. Consequently the Eagon-Northcott complex
associated to $\phi$
\[
    \sF^\bullet_\phi(d) :  \dots \to \sF^{-d-1} \boxtimes \sO(-d-1) \to \sF^{-d}\boxtimes \sO(-d)
\] 
with terms
\[
      \sF^{-i} = \Lambda^i \sG \otimes \Lambda^d \sH 
                 \otimes S_{i-d} \sH 
\]
gives a resolution
\[
       \sF^\bullet_\phi(d)[-d] \to I_X.
\]

\begin{lem}
$\sF^\bullet_\phi(d)$ is a family of (minimal degree) linear strands.
\end{lem}

\begin{proof}
Over each fiber $\sF^\bullet_\phi(d)$ restricts to the Eagon-Northcott complex that
resolves $I_{X_{|D|}}$.
\end{proof}

\begin{prop} \xlabel{p-psecd} 
 Let $d < n/2$ then we have the following identity of sheaves
\[
           \pi_* \sO_X = \sOSecd   
           \qquad \text{and} \qquad
           \pi_* I_X = I_\Secd.
\]
In particular $\cup X = \Secd$ as schemes.
\end{prop}

\begin{proof}
Set theoretically $\Sec_d = \cup X = \pi(X)$ since every $d$-secant of $E$ cuts out a degree $d$ 
divisor $D$ on $E$ and hence occurs in $X_{|D|}$. Since $\Secd$ is normal by
Corrolary \ref{c-secdnormal} and the restriction $\pi|_X \colon X \to \Sec_d$ is birational it follows that
$\pi_* \sO_X = \sOSecd$ (cf. \cite[proof of Corollary III.11.4]{Ha}). We apply $\pi_*$ to
the exact sequence
\[
         0 \to I_X \to \sO_{\Jacd \times \Pspace} \to \sO_X \to 0
\]
and obtain a long exact sequence
\xycenter{
         0 \ar[r]
         & I_{\cup X} \ar[r]
         & \sO \ar[rr] \ar[rd]
         & 
         & \sOSecd \ar[r]
         & R^1 \pi_* I_X \ar[r]
         & \dots
         \\
         &&& \sO_{\cup X} \ar[rd] \ar[ru]
         \\
         && 0 \ar[ru]
         && 0
         }
The fact that we have an injective morphism $\sO_{\cup X} \hookrightarrow \sOSecd$
shows that $\cup X$ is reduced and therefore $\cup X = \Sec_d$ as schemes.
 \end{proof}

\begin{prop} Let $d < n/2$ then we have the following identity of sheaves
\[
           R \sO_X = \sOSecdminus   
           \qquad \text{and} \qquad
           R I_X = I_\Secdminus.
\]
In particular $\cap X = \Secdminus$ as schemes.
\end{prop}

\begin{proof}
We shall first show that set theoretically $\cap X = \Secdminus$.
For this note that every $(d-1)$-secant $\PZ^{d-2}$ 
occurs in every $X_{|D|}$, we therefore
have $\Sec_{d-1} \subset \cap X$. Now if $n > 2d$ two $d$-secant 
$\PZ^{d-1}$'s can intersect in a $\PZ^{i-1}$ if and only if this
is an $i$-secant of $E$. Otherwise they would span a $\PZ^{2d-1-i}$
containing more than $2d-i$ points which is impossible by
Lemma \ref{l-normalspan}. Consequently we have $\Sec_{d-1} = \cap X$
as sets. 

For an arbitrary divisor $D$ of degree $d$ we consider the diagram
\xycenter{
            \tildeSecdminus \ar@{ (->}[r] \ar[dr]
            & \widehat{\Secdminus} \ar@{ (->}[r] \ar[d]
            & \widehat{\cap X} \ar@{ (->}[r] \ar[d]
            & \PE \ar[r]^{\tau} \ar[d]
            & \PZ^{d-1}
            \\
            &\Secdminus \ar@{ (->}[r]
            & \cap X \ar@{ (->}[r]
            & X_{|D|}
     }
where all squares are cartesian.
By Proposition \ref{p-canonicalclass} and Corrolary \ref{c-anticanonical} the preimage $\tildeSecdminus$ is a divisor on $\PE$ of class $dH - (n-2d)R$.

The intersection $\cap X$ is
cut out by hypersurfaces $Y_1, \dots Y_m$ of degree $d$ in $\Pspace$. Pulling these back to 
$\PE$ gives divisors $\widehat{Y_i}$ of class $dH$ containing $\widehat{\cap X}$
and hence also $\tildeSecdminus$. This shows that   
\[
     \widehat{Y_i} = \tildeSecdminus + \widehat{Z_i}
\]
with $\widehat{Z_i}$ a divisor of class $(n-2d)R$. In particular 
\[
      \widehat{Z_i} = \tau^{-1} Z_i
\]
where $Z_i$ is a hypersurface of degree $(n-2d)$ in $\PZ^{d-1}$. Altogether we find
\[
  \widehat{\cap X} = \tildeSecdminus \cup \tau^{-1} (\cap Z_i)
\]
as schemes. We claim that $\cap Z_i = \emptyset$. Otherwise set theoretically 
$\widehat{\cap X}$ and with it $\tildeSecdminus$ contains at least one fiber of $\tau$. The
image of this fiber is a $d$-secant $\PZ^{d-1}$ of $E$ contained in $\Sec_{d-1}$.
This is a contradiction since $\Sec_{d-1}$ contains no $d$-secant $\PZ^{d-1}$ of $E$ by Lemma
\ref{l-normalspan} in the appendix.

This shows $\tildeSecdminus = \widehat{\Secdminus} = \widehat{\cap X}$ and therefore also
$\Secd = \cap X$ as schemes. Consequently
\[
   R\sO_X = R^1\pi_*(\sO_X \otimes \omega_E) = \sO_{\cap X} = \sOSecdminus.
\]
As in the last proof we apply $\pi_*$ to
the exact sequence
\[
         0 \to I_X \to \sO_{\Jacd \times \Pspace} \to \sO_X \to 0
\]
to obtain the long exact sequence
\[
         0 \to I_\Secd \to \sO \xrightarrow{\alpha} \sOSecd \to R I_X \to \sO \to \sO_\Secdminus \to 0.
\]
Since $\alpha$ is surjective by the Proposition \ref{p-psecd} above this gives $R I_X = I_\Secdminus$.

\end{proof}

For future use we calculate the degree of the various vector bundles above.

\begin{lem}
Let $\sG$ and $\sH$ be as above. Then $\deg \sG = \deg \sH = -1$.
\end{lem}

\begin{proof}
We first have to give an explicit description of the Poincar\'e line bundle $\sP$. We denote
by $B \in \Jacd$ the point corresponding to the line bundle $\sOE(dA)$. By $\Delta$ we denote
the following divisor on $E \times \Jacd$
\[
      \Delta := \{ (P,\sL) \suchthat \sL = \sOE\bigl( (d-1) A + P \bigr) \}.
\]
One immediately checks that
\[
    \sP = \sOEJacd((d-1)p^*A + \Delta - q^* B).
\]
Since $\Jacd$ is an elliptic curve Riemann Roch gives $\chi(\sH) = \deg \sH$. By the Leray spectral sequence $\chi(\sH) = \chi(\sP)$ since $R^1 q_* \sP = 0$. Using Riemann Roch on the surface
$E \times \Jacd$ we obtain
\[
      \chi(\sP) = \frac{1}{2} ((d-1)p^*A + \Delta - q^* B)^2 = -1.
\]
The same reasoning works for
\[
     \sG = q_*(p^* \sOE(nA) \otimes \sP^{-1} ) = \sOEJacd((n-d+1)p^*A - \Delta + q^* B)
\]
and an analogue calculation gives
\[
       \deg \sG = \chi(\sG) = -1.
\] 
\end{proof}

Recall the following identities of exponential Chern characters on a smooth curve.

\begin{lem} \xlabel{l-Chernidentities}
Let $\sF$ be a rank $r$ bundle on a curve with exponential Chern character
\[
    \ch(\sF) = r + c_1(\sF).
\]
Then
\begin{align*}
    \ch(\Lambda^k \sF) &= {r \choose k} + {r-1 \choose k-1} c_1(\sF) \\
    \ch(S_k \sF)              &= {r + k - 1 \choose r-1} + {r+k-1 \choose r} c_1(\sF). \\
\end{align*}
\end{lem}

\begin{proof}
Well known, follows from standard calculation of chern classes on curves.
\end{proof}

\begin{lem} \xlabel{l-degFi}
Let $\sG$, $\sH$ and 
\[
      \sF^{-i} = \Lambda^i \sG \otimes \Lambda^d \sH 
                 \otimes S_{i-d} \sH 
\]
be as above. Then
\[
 \deg{\sF^{-i}}=-{n-d \choose i}{i \choose d} \frac{n}{n-d}.
\]
\end{lem}

\begin{proof}
Since $\sG$ and $\sH$ have ranks $(n-d)$ and $d$ respectively,
 the formulae of Lemma \ref{l-Chernidentities} give
\begin{align*}
     \ch(\Lambda^i \sG) &= {n-d \choose i} + {n-d-1 \choose i-1} c_1(\sG)\\
     \ch(\Lambda^d \sH) &= 1 + c_1(\sH) \\
     \ch(S_{i-d} \sH) &= {i-1 \choose d-1} + {i-1 \choose d} c_1(\sH). 
\end{align*}
By the multiplicativity of the exponential Chern character with respect to tensor products, we obtain
\begin{align*}
     \ch(\Lambda^d \sH \otimes S_{i-d} \sH) 
     &= {i-1 \choose d-1} + {i-1 \choose d} c_1(\sH) + {i-1 \choose d-1} c_1(\sH) \\
     &=  {i-1 \choose d-1} + {i \choose d} c_1(\sH) 
\end{align*}
and
\begin{align*}
\ch (\sF^{-i}) &= \ch(\Lambda^i \sG \otimes \Lambda^d \sH \otimes S_{i-d} \sH) \\ 
                      & =  {n-d \choose i} {i-1 \choose d-1} 
                              + {n-d \choose i}{i \choose d} c_1(\sH) 
                              + {n-d-1 \choose i-1}{i-1 \choose d-1} c_1(\sG) \\
                      & =  {n-d \choose i} {i-1 \choose d-1} 
                              + {n-d \choose i}{i \choose d}
                                  \left( 
                                      c_1(\sH) + \frac{i}{n-d}\frac{d}{i} c_1(\sG)
                                 \right) \\
                       & =  {n-d \choose i} {i-1 \choose d-1} 
                              + {n-d \choose i}{i \choose d}
                                  \left( 
                                      c_1(\sH) + \frac{d}{n-d} c_1(\sG)
                                 \right). 
\end{align*}
Since $\deg(\sH)=\deg(\sG)=-1$ this proves the lemma.
\end{proof}

\section{the geometric syzygy conjecture for $\Sec_{d-1}$} \xlabel{s-geosyz}
\nosubsections

\newcommand{\PG}{{\PZ(\sG)}}
\newcommand{\sOPG}{\sO_\PG}
\newcommand{\sErel}{{\sE_{rel}}}
\newcommand{\PErel}{{\PZ(\sErel)}}
\newcommand{\sOPErel}{\sO_\PErel}
\newcommand{\rhorel}{{\rho_{rel}}}

We have so far constructed the minimal free resolution of $\Sec_{d-1} := \Sec_{d-1}$
as a mapping-cone of the minimal free resolution of $X_{|D|}$ and its dual. 
This exhibits roughly half of the syzygies of $\Sec_{d-1} $ as geometric,
i.e. as coming from a scroll $X_{|D|}$ containing
$\Sec_{d-1}$. Now this construction can be repeated for all line bundles
in $\Jacd$ and we obtain different geometric syzygies each time. In this 
section we will show that the geometric syzygies obtained by this
procedure will generate the whole minimal free resolution of $\Sec_{d-1}$.

In the last section we construced a family of linear strands  
\[
     \sF^\bullet_\phi(d)[-d] \to I_X
\]
with $\cap X =\Sec_{d-1}$ and $\cup X = \Sec_d$.
By the procedure of Section \ref{s-intersection} we will obtain a complex 
$R\sF^\bullet$ that maps to the linear strand of $\cap X = \Sec_{d-1}$ 
and whose
image is a part of the linear strand that is generated by geometric
syzygies. We will prove that this map is surjective. The main ingredient
of this proof is the fact
that $\pi_*\sF^\bullet$ vanishes since $\cup X = \Sec_d$ does not
lie on any hypersurface of degree $d$.

We start by repeating our construction of the minimal free resolution of
$\Sec_{d-1}$ for all $X_{|D|}$ at once. The relative version of 
Corollary \ref{l-idealomega}
is

\begin{prop} \xlabel{p-relativeidealomega}
If $\omega_{X}$ is the dualizing sheaf of $X$ then
\[
    I_{\pi^{-1}(\Secdminus)/X} \cong \omega_X.
\]
\end{prop}

\begin{proof}
Consider the natural map
\[
     \sG \otimes \sH \to H^0(nA) \otimes \sOE
\]
obtained by fiberwise multiplication of sections. Projectivizing $\sG$ gives a short exact
sequence
\[
      0 \to  \sH \otimes \sOPG(-1) \to H^0(nA) \otimes \sOPG \to \sErel \to 0.
\]
The cokernel $\sErel$ is a relative version of the bundle $\sE$ considered in Section \ref{s-mfr}. 
Notice that $\PZ(\sG)$ is just the symmetric product $S^{d-1} E$. The projectivization $\PZ(\sErel)$ has appeared as the incidence varitey $I$ in the proof of \ref{p-secdminussmooth}.

We have a diagram
\xycenter{
                 & \PErel  \ar[d]^\rhorel \ar[dr]^{\sigma'} \ar[ddl]_{\pi'}
                 \\
                 & X \ar[r]^{\sigma} \ar[dl]^\pi
                 & \Jacd .
                 \\
                 \Pspace
                }
Recall that
\[
     \Secdminus \times \Jacd = \pi^{-1} (\Secdminus) .
\]
We set
\[
      \widetilde{\Secdminus \times \Jacd}  := (\pi')^{-1}(\Secdminus).
\]
Notice that over a point $|D| \in \Jacd$ this is just the preimage
\xycenter{
                  \widetilde{\Secdminus}  \ar[d] \ar@{ (->}[r]
                  & \PE   \ar[d]^\rho \ar@{ (->}[r]
                  & \PErel  \ar[d]^\rhorel
                  \\
                  \Secdminus \ar@{ (->}[r]
                  & X_{|D|} \ar@{ (->}[r]
                  & X
                  }
constructed in Section \ref{s-mfr}.

It follows from Corollary \ref{c-anticanonical} that
\[
     \sOPErel(-\widetilde{\Secdminus \times \Jacd}) 
     \cong 
     \omega_{\PErel} \otimes (\sigma')^* L
\] 
for some linebundle $L$ on $\Jacd$. Adjunction gives
\[
      \omega_{\widetilde{\Secdminus \times \Jacd} }
      \cong  (\sigma')^* L|_{\widetilde{\Secdminus \times \Jacd} }.
\]                
Pushing this forward with $\rhorel$ we obtain
\[
      \omega_{\Secdminus \times \Jacd }
      \cong  \sigma^* L|_{\Secdminus \times \Jacd}.
\]                
since $\Secdminus \times \Jacd$ is normal by Corollary $\ref{c-secdnormal}$.
Now $\omega_{\Secdminus \times \Jacd }$ is trivial by Corollary \ref{c-secdgorenstein} and pushing
down via $\sigma$ gives
\[
        \sOJacd \cong L.
\]
The same argument as in Lemma \ref{l-idealomega} then proves our claim.
\end{proof}

\begin{prop}
There is a short exact sequence
\[
      0 \to \sO_{\Jacd \times \Pspace}(-n) \to (\sF^\bullet_\phi(d))^*[-d] \to \omega_{X} \to 0
\]
with $(\sF^\bullet_\phi(d))^*=\sHom(\sF^\bullet_\phi(d), \omega_{\Jacd \times \Pspace})[n]$.
\end{prop}

\begin{proof}
Applying $\sHom(-,\omega_{\Pspace \times \Jacd})$ to the exact sequence
\[
           0 \to \sF^\bullet[d] \to \sO_{\Jac_d \times \Pspace} \to \sO_{X}
\]
yields the desired result.
\end{proof}

\begin{cor} \xlabel{c-dualkoszul}
Let $d \le i \le n-d$ then
\[
(\sF^{-n+i})^* = \sigma_* (
               \omega_{X} \otimes \pi^* \Omega^{i-d}(i)
               ).
\]
\end{cor}

\begin{proof}
Since 
\[
      0 \to \sO(-n) \to (\sF^\bullet_\phi(d))^*[-d] \to \omega_{X} \to 0
\]
is an exact resolution of $\omega_X$ and 
\[
               R^q \sigma_* (\pi^*\sO(-n) \otimes \pi^* \Omega^{i-d}(i) )
              = H^q (\Omega^{i-d}(i-n)) \otimes \sOJacd = 0
\]
in the range $d\le i \le n-d$, the same spectral sequence argument 
as in \ref{l-relativekoszul} applies.
\end{proof}

Each choice of $|D| \in \Jacd$ exhibits the minimal free resolution
of $\Sec_{d-1}$ as a mapping cone by Theorem \ref{t-mfr}. 
The information how these mapping cones fit together is
encoded in the short exact sequences of the  following proposition.

\begin{prop}
Let $d \le i \le n-d$. Then
there are short exact sequences
\[
      0 \to \sF^{-i} \to F^{-i}_{\Sec_{d-1}} \otimes \sOJacd \to (\sF^{-n+i})^* \to 0
\]
on $\Jac_{d}$.
\end{prop}

\begin{proof}
Consider the short exact sequence
\[
       0 \to I_{X/ \Jacd \times \Pspace}
       \to I_{\Jac_d \times \Sec_{d-1}/ \Jacd \times \Pspace}
       \to I_{\Jac_d \times \Sec_{d-1}/X} 
       \to 0.
\]
We tensor this sequence with $\pi^*\Omega^{i-d}(i)$ and
apply $\sigma_*$. By Theorem \ref{t-mfr} we know that the postulation
of $\Secdminus$ is $d$. Lemma \ref{l-relativekoszul} then tells us that
\[
  \sigma_* (I_X \otimes \pi^* \Omega^{i-d}(i)) = \sF^{-i}
 \quad \text{and} \quad
   R^1 \sigma_* I_X \otimes \pi^* \Omega^{i-d}(i) = 0.
\]
Lemma \ref{l-koszul} gives us
\[
     \sigma_* (I_{\Jac_d \times \Sec_{d-1}/ \Jacd \times \Pspace} \otimes \pi^* \Omega^{i-d}(i)) 
     = F^{-i}_\Secdminus \otimes \sOJacd.
\]
By Proposition \ref{p-relativeidealomega}
\[
     I_{\Jac_d \times \Sec_{d-1}/ X} \cong \omega_X
\]
and now Corollary \ref{c-dualkoszul} shows
\[
      \sigma_* (\omega_X \otimes \pi^* \Omega^{i-d}(i))
      =(\sF^{-n+i})^* 
\]
This proves the result.
\end{proof}

We shall also need the following.

\begin{prop} \xlabel{p-vanishing}
Let $d \le i \le n-d$, then
$H^0(\Jac_d,\sF^{-i})=0$.
\end{prop}

\begin{proof}
By Theorem \ref{t-mfr} we know, that the postulation of $\Secd$
is $d+1$. Therefore by Proposition \ref{p-cupstrand}
we have
\[
  \pi_* \sF_\phi^\bullet = 0,
\] 
in particular
\[
      (\pi_* \sF_\phi^\bullet)_i = H^0(\sF^{-i}) \otimes \sO(-i-d) = 0.
\]
\end{proof}

\begin{thm}
The linear syzygies of $\Sec_{d-1}$ are spanned by scrollar syzygies.
\end{thm}

\begin{proof}
We recall the map $\psi \colon R\sF^\bullet(d) \to F^\bullet_{\Sec_{d-1}}$ from Proposition \ref{p-psi}
whose image is the part of the minimal free resolution of $\Sec_{d-1}$ generated by scrollar
syzygies. 
We have to show that $\psi$ is surjective. For this it is enough to check that in each step the natural map
\[
    \psi_{-i} \colon H^1(\sF^{-i}) \to F^{-i}_{\Sec_{d-1}}
\]
is surjective. For this consider the short exact sequence
\[
      0 \to \sF^{-i} \to F^{-i}_{\Sec_{d-1}}  \otimes \sOJacd \to (\sF^{-n+i})^* \to 0.
\]
The associated long exact sequence shows that the claim is true,
if $H^1((\sF^{-n+i})^*)=0$. But this follows from Proposition \ref{p-vanishing}
since
\[
    H^1((\sF^{-n+i})^*)=H^0(\sF^{-n+i})
\]
by Serre duality on the elliptic curve $\Jac_d(E)$.
\end{proof}

We now ask wether the natural map 
$\psi_{-i} \colon H^1(\sF^{-i}) \to F^{-i}_{\Sec_{d-1}}$ is an isomorphism. 
Since we already
know the dimension of $F^{-i}_{\Sec_{d-1}}$ by Theorem \ref{t-mfr}, we only need
to calculate $H^1(\sF^{-i}) = - \deg \sF^{-i}$ since $H^0(\sF^{-i})$ vanishes
and $E$ is elliptic.
\begin{cor} \xlabel{c-isom}
$\psi \colon H^1(\sF^{-i}) \to F^{-i}_{\Sec_d}$ is an isomorphism 
for $i=d$ and $i=n-d$.
\end{cor}

\begin{proof}
Since $\psi$ is surjective, we only have to compare the dimensions
of the vector spaces. By Lemma \ref{l-degFi} we obtain for $i=d$
\[
      H^1(\sF^{-d}) = - \deg \sF^{-d} = {n-d \choose d}{d \choose d} \frac{n}{n-d}
\]
since $H^0(\sF^{-d})=0$. On the other hand Corollary \ref{c-secdbetti} yields
\begin{align*}
    \dim F_{\Secdminus}^{-d} 
     &= { n-d \choose d}{d-1 \choose d-d}
       +{ n-d \choose n-d}{n-d-1 \choose n-d-d} \\
     &=  { n-d \choose d}\frac{n}{n-d} \\
\end{align*}
as claimed. For $i=n-d$ we obtain the same numbers.
\end{proof}

\begin{rem} \xlabel{r-notlinearstrand}
For other steps in the resolution $\psi$ is never an isomorphism.
For example if $n=6$ the the Betti  diagram of  $R\sF^\bullet$ and $E$ are
respectively
\[
\begin{matrix}
     \begin{matrix}
       - & - & 9 & 18 & 9 & - & - \\
     \end{matrix}

     &
     \quad\text{and}\quad
     &

     \begin{matrix}
       - & - & -  & - & 1 \\
       - & 9 & 16 & 9 & - \\
       1 & - & -  & - & -
     \end{matrix}
\end{matrix}
\]
The middle term of $R\sF^\bullet$ is to large.
In the next section we will explain this difference.
\end{rem}

\section{The kernel of $\psi$} \xlabel{s-kernel}

As we have seen, the linear complex $R \sF_\phi^\bullet$ is
not the linear strand of $\Sec_{d-1}$. In this section we
will describe the difference. As it turns out, this difference
is the linear strand of $\Sec_{d}$.

\begin{prop}
There exist an exact sequence
\[
         0 \to F_{\Sec_d}^\bullet(d+1) 
           \to R \sF_\phi^\bullet
           \xrightarrow{\psi} F_{\Sec_{d-1}}^\bullet(d)
           \to 0
\]
\end{prop}

\begin{proof}
Let $\ker^\bullet$ be the kernel of the surjective morphism 
\[
    \psi \colon R \sF_\phi^\bullet(d) \to F_{\Sec_{d-1}}^\bullet(d).
\]
Since $\psi_{-d}$ and $\psi_{-n+d}$ are isomorphisms by 
Corollary $\ref{c-isom}$,
we have $\ker^{-d}=\ker^{-n+d}=0$, and a commutative
diagram
\begin{center}
\mbox{
\xymatrix{      
             0  \ar@{ (->}[r] 
           & H^1(\sF_\phi^{-d})\otimes\sO(-d) 
             \ar@{->>}[r]^-{\psi_{-d}} 
           & F_{\Sec_{d-1}}^{-d}\otimes\sO(-d) 
           \\
             \ker^{-(d+1)}   \ar@{ (->}[r] \ar[u]
           & H^1(\sF_\phi^{-d-1})\otimes\sO(-d-1)  
            \ar@{->>}[r]^-{\psi_{-d-1}} \ar[u]
           & F_{\Sec_{d-1}}^{-d-1}\otimes\sO(-d-1)  \ar[u]
           \\
             \dots  \ar[u]
           & \dots  \ar[u]
           & \dots  \ar[u]
           \\
             \ker^{-n+(d+1)}   \ar@{ (->}[r] \ar[u]
           & H^1(\sF_\phi^{n-d-3})\otimes\sO(-n+d+1)  
             \ar@{->>}[r] 
           & F_{\Sec_{d-1}}^{-n+d+1}\otimes\sO(-n+d+1)  \ar[u]
           \\
             0  \ar@{ (->}[r] \ar[u]
           & H^1(\sF_\phi^{-n+d})\otimes\sO(-n+d)  
             \ar@{->>}[r]^-{\psi_{-n+d}} \ar[u]
           & F_{\Sec_{d-1}}^{-n+d}\otimes\sO(-n+d)  \ar[u]
          }
     }
\end{center}
 To show that $\ker^\bullet$ is the degree $d+1$ linear
strand of $\Sec_d$, we calculate the homology of the vertical complexes 
above.

For the middle column $R\sF_\phi^\bullet(d)$ 
we use the hypercohomology with respect to $\pi$.
Since $\sF^\bullet_\phi(d)$ has cohomology $I_X$ concentrated in step $-d$,
the second hypercohomology spectral sequence
\[
      {}^{II} E_{2}^{pq} = R^p \pi_* (H^q(\sF_\phi^\bullet(d))) 
                \Rightarrow \HZ^{p+q}_{\pi}(\sF_\phi^\bullet(d))
\]
yields
\[
    \HZ^{p+q}_{\pi}(\sF_\phi^\bullet(d)) 
    = \left\{
    \begin{array}{cc}
          \pi_* I_X = I_\Secd & \text{if $p+q=-d$} \\
          R^1 \pi_* I_X = I_\Secdminus & \text{if $p+q = -d+1$} \\
          0 & \text{otherwise} .
    \end{array}
    \right.
\]
Consider now the first hypercohomology spectral sequence
\[
      {}^I E_{2}^{pq} = H^p(R^q \pi (\sF_\phi^\bullet(d) )) 
                \Rightarrow \HZ_\pi^{p+q}(\sF_\phi^\bullet(d)).
\]
Since
\[
    R^q\pi_* (\sF^{-i} \boxtimes \sO(-i)) 
    = H^q(\sF^{-i}) \otimes \sO(-i)
    = \left\{
    \begin{array}{cc} 
           H^1(\sF^{-i})\otimes \sO(-i) & \text{if $q=1$} \\
                    0                                     & \text{otherwise}
    \end{array}
    \right.
\]
this spectral sequence degenerated in level 2 and we obtain
\[
        \HZ_\pi^{p+1}(\sF^\bullet_\phi(d)) = H^p(R^1\pi_* \sF^\bullet_\phi(d)).
\]
Comparing this with the result from the second hypercohomology sequence
shows that $R\sF^\bullet_\phi(d)$ has cohomology $I_\Secd$ and $I_\Secdminus$
concentrated in steps $-d-1$ and $-d$.

For the right column we notice, that the minimal free resolution
of $I_{\Sec_{d-1}}$ is
\[
      0 \to \sO(-n) \to F_{\Sec_{d-1}E}^\bullet(d) \to I_{\Sec_{d-1}} \to 0
\]
and therefore $F_\Secdminus^\bullet(d)$ has cohomology
$I_{\Sec_{d-1}}$ in step $-d$ and 
$\sO(-n)$ in step $-n+d$.

All together we obtain the long exact cohomology sequence

\begin{center}
\mbox{
\xymatrix{      
           & 0 \ar[r]
           & I_{\Sec_{d-1}}  \ar[r]^{\cong}
           & I_{\Sec_{d-1}}  
           &
           \\
           & H^{-d-1}(\ker^\bullet) \ar[r]^-{\cong}
           & I_{\Sec_d}  \ar[r]
           & 0 \ar `r[u] `[lll] `[ulll]  `[ull] [ull]
           &
           \\
           & H^{-d-2}(\ker^\bullet) 
           & 0 \ar[r] 
           & 0  \ar `r[u] `[lll] `[ulll]  `[ull] [ull]
           &
           \\
           & \dots  
           & \dots  
           & \dots  
           &  
           \\
           & H^{-n+d+1}(\ker^\bullet)  \ar[r] 
           & 0 \ar[r] 
           & 0  
           & 
           \\
           & 0  \ar[r] 
           & 0 \ar[r] 
           & \sO(-n)  \ar `r[u] `[lll] `[ulll]_-\cong  `[ull] [ull] 
           &  
          }
     }
\end{center}

This shows that $\sO(-n) \to \ker^\bullet[d+1] \to I_{\Sec_d}$ 
is exact and
therefore $$F_{\Sec_d}^\bullet(d+1) = \ker^\bullet$$.
\end{proof}
 
\begin{example}
For $n=6$ we can write this suggestively as
\[
  \begin{matrix}
      - & - & -  & - & 1 \\
      - & 9 & 16 & 9 & - \\
      1 & - & -  & - & - \\
  \end{matrix}
  \quad
  +
  \quad
  \begin{matrix}
      - & - & 1 \\  
      - & - & - \\
      - & 2 & -  \\
      - & - & -  \\
      1 & - & - \\
  \end{matrix}
  \quad = \quad
  \begin{matrix}
      - & - & -  & \not1 & - \\
      - & - & -  & - & \not1 \\
      - & 9 & 18 & 9 & - \\
  \not1 & - & -  & - & - \\
      - & \not1 & -  & - & - \\
  \end{matrix}
\]
(In a Betti diagram all entries on a diagonal have the same twist,
therefore one can sometimes ``cancel diagonally''.)
\end{example}

\begin{example}
For $n=7$ we write similarly
\[
  \begin{matrix}
      - & -  & -  & -  & -  & 1 \\
      - & 14 & 35 & 35 & 14 & - \\
      1 & -  & -  & -  & -  & - \\
  \end{matrix}
  \quad 
  +
  \quad
  \begin{matrix}
      - & - & - & 1 \\  
      - & - & - & - \\
      - & 7 & 7 & -  \\
      - & - & - & - \\
      1 & - & - & - \\
  \end{matrix}
  \quad
  =
  \quad
  \begin{matrix}
      - & -  & -  & -  & \not1 & - \\
      - & -  & -  & -  & -  & \not1 \\
      - & 14 & 42 & 42 & 14 & - \\
      \not1 & -  & -  & -  & -  & - \\
      - & \not1  & -  & -  & -  & - \\
  \end{matrix}
\]
\end{example}

\section{The geometric Syzygy Conjecture for bielliptic canonical curves} \xlabel{s-bielliptic}
\nosubsections

As an application of our study of elliptic normal curves, we deduce
the geometric syzygy conjecture for bielliptic canonical curves. This will follow
directly from the classic fact, that a bielliptic canonical curve $C$
is the transversal intersection of a quadric $Q$ in $\PZ^{g-1}$ with a cone $\tilde{E}$
over an elliptic normal curve $E$ of degree $g-1$. From this we obtain the
minimal free resolution of $\sO_C$ as the total complex of
\[
         \bigl[F^\bullet_Q[1] \to \sO_{\PZ^{g-1}}\bigr] 
         \otimes \bigl[F^\bullet_{\tilde{E}}[1] \to \sO_{\PZ^{g-1}}\bigr]
\]
where $F^\bullet_Q$ and $F^\bullet_{\tilde{E}}$ are the minimal free 
resolutions of $I_Q$ and $I_{\tilde{E}}$. In particular, all higher syzygies
of $C$ come from syzygies of $\tilde{E}$, which in turn come from geometric
syzygies.

Let us start by recalling the classical facts.

\begin{prop}
Let $C \subset \PZ^{g-1}$ be a bielliptic canonical curve of
genus $g$ and Clifford index $c \ge 2$. Let $\iota \colon C \xrightarrow{2:1} E$ be the elliptic involution.
Then
\[ 
      \tilde{E} = \bigcup_{e \in E} \spans{\iota^{-1}(e)}
\]
is an elliptic cone with vertex $p \not\in C$. The projection of $\tilde{E}$
from $p$ is an elliptic normal curve $E \subset \PZ^{g-2}$. Moreover there
exists a (not uniquely defined) quadric $Q \subset \PZ^{g-1}$ such that
\[
    C = \tilde{E} \cap Q
\]
and this intersection is transversal.
\end{prop}

\begin{proof}
The fact that all lines $\spans{\iota^{-1} e}$ contain a common point
$p \not\in C$ is \cite[p. 269, Exercise E-1]{ACGH}. Projecting
from this point gives an elliptic curve $E \subset \PZ^{g-2}$. The
degree of $E$ is $g-1$ since the projection also induces
the $2:1$ map from the degree $2g-2$ curve $C$ to $E$. Consequently
$E$ is an elliptic normal curve, and $I_{\tilde{E}}$ is the
pullback of $I_E$ via the projection from $p$. From Theorem \ref{t-mfr} respectively 
Corollary \ref{c-secdbetti} applied to $E=\Sec_1$ we know
that $I_E$ is generated by 
\[
   { g-3 \choose 2}\frac{g-1}{g-3}= { g-2 \choose 2} - 1
\] 
quadrics. On the other hand by Petri's Theorem \cite[p.131]{ACGH} 
the ideal of every non-trigonal canonical 
curve without a $g^2_5$ is generated by ${g-2 \choose 2}$ quadrics. 
Consequently there exist a quadric $Q$ that contains $C$ but not $\tilde{E}$.
The intersection $\tilde{E} \cap Q$ is then transversal and yields $C$.
\end{proof}

\begin{cor}
The minimal free resolution of $C \subset \PZ^{g-1}$ is the total complex
of
\[ 
     [\sO_{\PZ^{g-1}}(-2) \xrightarrow{Q} \sO_{\PZ^{g-1}}] \otimes 
     [\sO_{\PZ^{g-1}}(-g-1) \to \iota^* F^\bullet_E(2)[-1] \to \sO_{\PZ^{g-1}}]
\]
where $F^\bullet_E(2)$ is the linear strand of $E \subset \PZ^{g-2}$.
\end{cor}

\begin{proof}
Since $\tilde{E}$ is a cone over $E$, the pullback of
the minimal free resolution of $E$
\[
     \sO_{\PZ^{g-2}}(-g+1) \to F^\bullet_E(2)[-1] \to \sO_{\PZ^{g-2}}
\]
via the projection from the vertex $p$ gives the minimal free resolution of $\tilde{E}$. 
Now the intersection $C = Q \cap \tilde{E}$ is transversal,
so the total complex above gives a free resolution of $C$. Since
$C \subset \tilde{E}$ and both have linear the same postulation,
$F^\bullet_{\tilde E}(2)$ must be a subcomplex of the minimal free
resolution of $C$. Therefore no canceling in the above total complex
can occur, i.e. the total complex is a minimal free resolution of $\sO_C$.
\end{proof}

\begin{example}
For $g=7$ we obtain the following Betti  diagrams
\[
     \begin{matrix}
       \begin{matrix}
            - & 1 \\
            1 & - 
       \end{matrix}
      &
       \otimes
      & 
        \begin{matrix}
          - & - & -  & - & 1 \\
          - & 9 & 16 & 9 & - \\
          1 & - & -  & - & - \\
        \end{matrix}
      &
       =
      &
        \begin{matrix}
          - & -  & -  & -  &  - & 1 \\
          - & -  & 9  & 16 & 10 & - \\
          - & 10 & 16 & 9  & -  & - \\
          1 & -  & -  & -  & -  & - \\
        \end{matrix}
      \\
      Q
      &&
      \tilde{E}
      &&
      C
   \end{matrix}
\]
\end{example}

\begin{cor}
Let $C \in \PZ^{g-1}$ be a bielliptic canonical curve of Clifford index $c \ge 2$.
Then the geometric syzygy conjecture is true for every step of the
minimal free resolution of $C$.
\end{cor}

\begin{proof}
By a theorem of Green \cite{GreenQuadrics} the rank $4$ quadrics in
the ideal of any canonical curve $C$ span the degree $2$ part of $I_C$.
This is the geometric syzygy conjecture in step $0$. By the preceding
proposition all other terms in the linear strand $F^\bullet_C(2)$ come
from the minimal free resolution of $E$ and are therefore spanned by geometric
syzygies.
\end{proof}

\begin{rem}
More precisely we have shown, that the higher linear 
syzygies of $C$ are generated by syzygies of scrolls constructed
from the $g^1_4$'s of $C$ that factor over a $g^1_2$ on $E$.
\end{rem}

\begin{rem}
The geometric syzygy conjecture is known for curves of Clifford index $1$, i.e.
trigonal curves \cite{Eusen} and curves which are isomorphic to plane quintics. 
For trigonal curves it is again enough to consider scrollar syzygies, but for plane quintics the scrollar syzygies no longer suffice and one has to take other geometric syzygies into account. See for example \cite[p. 87, Bemerkung 5.5.4]{HC}.
\end{rem}

\section{apendix}
\nosubsections

\begin{lem} \xlabel{l-cohomology}
The cohomology of $\Omega^i(i+j)$ is concentrated in 
one degree. More explicitly it is zero except for the following
cases:
\[
    \begin{array}{|c|c|}
    \hline
    \text{twist} & \text{non zero} \\
                 & \text{cohomology} \\
    \hline
     j<-n-1   & H^n \\
     j=-i     & H^i \\
     j>0      & H^0 \\
    \hline
    \end{array}
\]
Moreover $H^i(\Omega^i)=1$.
\end{lem}

\begin{proof} Use the Theorem of Bott. \end{proof}

\begin{lem} \xlabel{l-normalspan}
Let $E \subset \mathbb{P}^{n-1}$ be an elliptic normal curve
of degree $n$. Then every (possibly non-reduced) divisor $D$ of
degree $d \leq n-1$ spans a linear
subspace of dimension $d-1$
\end{lem}

\begin{proof}
Assume that this statement is false. Then given any $n-d$ points $P_1,\ldots,
P_{n-d}$ there exists a hyperplane $H$ which contains $D$ and the points
$P_1,\ldots, P_{n-d}$. On the other hand $D + P_1 + \ldots P_{n-d}$ is not
linearly equivalent to a hyperplane section if the $P_i$ are chosen
generically. This gives a contradiction.
\end{proof}



\end{document}